\documentclass[11pt]{article}
\usepackage[margin=1in]{geometry}  
\usepackage{graphicx}
\usepackage{subcaption}              
\usepackage{amsmath}               
\usepackage{amsfonts}              
\usepackage{amsthm}
\usepackage[square, comma, numbers, sort & compress]{natbib} 
\usepackage{hyperref}
\usepackage{mathtools}
\usepackage{tikz}
\usetikzlibrary{decorations.pathmorphing}
\usetikzlibrary{shapes}
\usetikzlibrary{arrows,decorations.pathreplacing}
\usepackage{circuitikz}
\usepackage{float}

\allowdisplaybreaks
\usepackage{bigints} 
\usepackage{resizegather}

\usepackage{epstopdf}
\newtheorem{Def}{Definition}
\newtheorem{Thm}{Theorem}

\newtheorem{Prop}{Proposition}
\newtheorem*{Pf}{Proof}
\newtheorem{Lem}{Lemma}
\newtheorem{Rk}{Remark}
\newtheorem{Cc}{Corollary}
\newtheorem{Not}{Notation}

\newtheorem{Asmpt}{Assumption}

\begin{document}

\title{A boundary feedback analysis for input-to-state-stabilisation of non-uniform linear hyperbolic systems of balance laws with additive disturbances}

\author{Mapundi K. Banda, Gediyon Weldegiyorgis}

\maketitle
\begin{abstract}
A boundary feedback stabilisation problem of non-uniform linear hyperbolic systems of balance laws with additive disturbance is discussed. A continuous and a corresponding discrete Lyapunov function is defined. Using an input-to-state-stability (ISS)  $ L^2- $Lyapunov function, the decay of solutions of linear systems of balance laws is proved. In the discrete framework, a first-order finite volume scheme is employed. In such cases, the decay rates can be explicitly derived. The main objective is to prove the Lyapunov stability for the $L^2$-norm for linear hyperbolic systems of balance laws with additive disturbance both analytically and numerically. Theoretical results are demonstrated by using numerical computations.
\end{abstract}


AMS:\\
  65Kxx, 49M25, 65L06


\section{Introduction}\label{sec:intro}

We consider the following $ k \times k $ non-uniform linear hyperbolic system of balance laws with additive disturbances (see Equation 2 in \cite{prieur2012iss} and Equation 1.1.10 in \cite{strikwerda2004finite}):
\begin{equation}\label{eq:LHSBLs}
\partial_t W + \Lambda (x) \partial_x W + \Pi(x) W = \Psi(x,t),\; x \in [0,l],\; t \in [0,+\infty),
\end{equation}
where $ W := W(x,t): [0,l]\times [0, +\infty) \rightarrow \mathbb{R}^k $ is a state vector. In addition $\Lambda(x) = \text{diag}\{{\Lambda^+}(x), -{\Lambda^-}(x) \} $, where $ {\Lambda^+}(x) \in \mathbb{R}_{+}^{m \times m}$ and $ {\Lambda^-}(x) \in \mathbb{R}_{+}^{(k-m) \times (k-m)}$, are non-zero differentiable diagonal matrices,  $ \Pi(x) \in \mathbb{R}^{k\times k} $  is a non-zero matrix and $ \Psi:=\Psi(x,t):[0,l]\times [0, +\infty) \rightarrow \mathbb{R}^k $ is a vector of disturbance functions. Corresponding to the positive and negative diagonal entries of $ \Lambda(x) $, the state vector $ W $ is specified by $ W^\top = [{W^+}, {W^-}]^\top $, where $ {W^+} \in \mathbb{R}^{m}$ and $ {W^-} \in \mathbb{R}^{k-m}$ and the disturbance function is also written as $ \Psi^\top = [{\Psi^+}, {\Psi^-}]^\top $, where $ {\Psi^+} \in \mathbb{R}^{m}$ and $ {\Psi^-} \in \mathbb{R}^{k-m}$. More clarity on the notation will be presented in Section \ref{sec:sec01}.

Equation \eqref{eq:LHSBLs} is supplemented by an initial condition which is set as
\begin{equation}\label{eq:ICLHSBLs}
W(x,0) = W_0(x), \; x \in (0,l),
\end{equation}
where $W_0 : (0,l) \rightarrow \mathbb{R}^k $ is of class $ C^1 $. On a finite spatial domain, the following linear feedback boundary conditions with no additive disturbance are prescribed: 
\begin{equation}\label{eq:BCsLHSBLs}
\begin{bmatrix} W^+(0,t) \\ W^-(l,t) \end{bmatrix} = K \begin{bmatrix} W^+(l,t) \\ W^-(0,t)  \end{bmatrix}, \; t > 0,
\end{equation}
where $ K \in \mathbb{R}^{k\times k} $ is a non-zero real matrix of the form 
$ K = \begin{bmatrix} 0 & K^- \\ K^+ & 0 \end{bmatrix} $, with $ {K^-} \in \mathbb{R}^{m \times (k-m)} $ and $ {K^+} \in \mathbb{R}^{(k-m) \times m} $, together with a zero-order initial boundary compatibility condition expressed as: 
\begin{equation}\label{eq:CCs01LHSBLs}
\begin{bmatrix} W^+(0,0) \\ W^-(l,0) \end{bmatrix} = K \begin{bmatrix} W^+(l,0) \\ W^-(0,0)  \end{bmatrix}. 
\end{equation}

The main purpose of this paper is to analyse numerical boundary feedback stability of non-uniform linear hyperbolic systems of balance laws with additive disturbances such as presented in equations \eqref{eq:LHSBLs} - \eqref{eq:CCs01LHSBLs} above. A numerical ISS $L^2- $Lyapunov function is constructed and used to investigate conditions for ISS. Mathematical proofs for decay rates for an upwind finite-volume scheme using an equivalent discrete Lyapunov function will be presented. The secondary purpose is to analyse the decay of a continuous ISS $L^2- $ Lyapunov function for the same hyperbolic systems. This serves as motivation to investigate the conditions under which in a numerical scheme decay of the discrete solutions can be achieved.  In addition, the decay of the discrete ISS $L^2- $Lyapunov function is also confirmed using numerical computations of linear hyperbolic systems and the Saint-Venant system.

Boundary feedback stabilisation of hyperbolic systems of balance laws, in general, has been an active research field,  see \cite{bastin2011boundary, krstic2008backstepping, coron2007strict, de2003boundary, dos2008boundary, diagne2012lyapunov, christofides1996feedback, bastin2008using, coron2015dissipative, gugat2014boundary, gugat2011existence}, for some references. However, boundary feedback stabilisation of linear hyperbolic systems of balance laws with additive disturbance is more recent \cite{prieur2012iss, gugat2018boundary}. For hyperbolic systems of balance laws, a strict $ L^2- $Lyapunov function is used to investigate conditions for exponential stability of such systems. In \cite{prieur2012iss} an ISS $L^2- $Lyapunov function is used to investigate conditions for ISS of time-varying linear hyperbolic system of balance laws with additive disturbance. In the current article, an ISS $L^2-$ Lyapunov function for non-uniform linear hyperbolic systems of balance laws is introduced in Section \ref{sec:sec01}. Therein a rigorous discussion of the decay of such a Lyapunov function is discussed.

It must be mentioned that in the field of dynamical systems control, input-to-state stability (ISS) is well known and it is used to analyse stability of nonlinear dynamical systems with additive disturbances (or external inputs) \cite{sanchez1999input, sontag1995characterizations, hespanha2008lyapunov}. For further study of ISS,  the reader is referred to \cite{sontag2008input}.    

Numerical boundary feedback stabilisation of hyperbolic systems of balance laws has become a developing research field \cite{banda2013numerical, dick2014stabilization, gottlich2017numerical, Banda2018, herty2016boundary, gersterdiscretized, gottlich2016electric}. In these studies, a discrete $ L^2- $Lyapunov function is constructed and used to investigate conditions for exponential stability of discretised hyperbolic systems. This is the main thrust of this paper. The numerical boundary feedback is discussed in Section \ref{sec:sec02} which contains the main results of this paper. Furthermore, the decay of the discrete ISS $ L^2- $Lyapunov function has been rigorously proved. The main contribution of this paper is a new numerical Lyapunov function and proof of its decay. In most cases such as \cite{prieur2012iss} numerical approaches are applied and their results are computationally demonstrated without a numerical analysis. This paper intends to fill that gap.

Exponential decay of the strict $ L^2- $Lyapunov function has been shown for some important physical problems such as gas flow through a pipeline \cite{gugat2011existence}, the transmission of electricity along a power line (defined by the telegraph equation) \cite{gugat2014boundary} and the shallow water flow along a channel with and without transportation of sediment \cite{bastin2011boundary, diagne2012lyapunov}. In this paper, the theoretical and numerical results are applied to a non-uniform linear system of balance laws with additive disturbances as well as to the well known Saint-Venant equations. The results in Section \ref{sec:sec03} demonstrate how the analysis can be applied and the numerical decay of the Lyapunov function can be observed. 

\section{Boundary feedback for ISS}\label{sec:sec01}
In this section the analytical boundary feedback results are presented and proved. Some necessary notation and definitions will be presented first and the main theorem of the section will be presented and proved. The section ends with a corollary which links the results herein with ISS for uniform linear balance laws. 
\begin{Not}
	Denote $ \mathbb{R}^{k} $, $ \mathbb{R}^{k \times k} $ and $ \mathbb{R}_{+}^{k \times k} $ as the set of $ k- $dimensional real vectors, $ k- $dimensional square real matrices and $ k- $dimensional square real matrices with positive entries, respectively. Denote $C^0$ and $C^1$ as the set of continuous and  continuously differentiable functions in $ \mathbb{R}^{k} $, respectively. For a given function $ f: [0,l] \rightarrow \mathbb{R}^{k} $, the $ L^2- $norm is defined as 
	\begin{equation*}
	\| f \|_{L^2} = \sqrt{\int_{0}^{l} |f(x)|^2dx}, 
	\end{equation*}
	where $ |\cdot| $ is the Euclidean norm in $ \mathbb{R}^{k} $. Moreover, $ L^2(0,l) $ is called the space of all measurable functions, $ f $, for which $ \| f \|_{L^2} < \infty $.
\end{Not}
In addition, the following assumptions are made:
\begin{Asmpt}\label{Asmpt:LHSBLs-Asmpt01}
For all $ x \in [0,l]$, and $ t \in [0,+\infty)$, assume 
\begin{enumerate}
	\item[(i)] The real diagonal matrix $ \Lambda $ is of class $ C^1 $.
\item[(ii)] The real matrix $ \Pi $ and the disturbance function $ \Psi $ are of class $ C^0 $.
	\item[(iii)] The $\displaystyle \sup_{s \in [0,t]}\left(|\Psi(x, s)|^2\right) $ is sufficiently small. 
\end{enumerate}
\end{Asmpt}

The existence and uniqueness of a solution to the system \eqref{eq:LHSBLs} with initial conditions \eqref{eq:ICLHSBLs}, boundary conditions \eqref{eq:BCsLHSBLs} and compatibility conditions \eqref{eq:CCs01LHSBLs} is proved in \cite{kmit2008classical}. The ISS of a steady-state $W\equiv 0 $ is defined as follows: 

\begin{Def}[ISS]
The steady-state $ W(x,t) \equiv 0 $ of the system \eqref{eq:LHSBLs} with the boundary conditions \eqref{eq:BCsLHSBLs} is ISS in $L^2-$norm with respect to disturbance function $ \Psi $ if there exist positive real constants $\eta > 0$, $ \xi > 0 $ and $ C > 0$ such that, for every initial condition $W_0(x) \in L^2((0,l);\mathbb{R}^k)$  satisfying the compatibility condition \eqref{eq:CCs01LHSBLs}, the  $L^2-$solution to the system \eqref{eq:LHSBLs} with initial condition \eqref{eq:ICLHSBLs} and boundary conditions \eqref{eq:BCsLHSBLs} satisfies
\begin{equation}\label{eq:ExponStabCondChap1}
	{\|W(\cdot, t) \|}_{L^2((0,l);\mathbb{R}^k)}^2 \leq C {e}^{-\eta t}{\|W_0 \|}_{L^2((0,l);\mathbb{R}^k)}^2 + \frac{C}{\eta \xi} \sup_{s \in [0,t]} \left(\int_{0}^{l}|\Psi(x,s)|^2dx\right),\;t \geq 0.
\end{equation}	
\end{Def}
\begin{Rk}
The extra term in the inequality \eqref{eq:ExponStabCondChap1} estimates the influence of the disturbance $ \Psi $ on the solution of the system \eqref{eq:LHSBLs} with the boundary conditions \eqref{eq:BCsLHSBLs}.
\end{Rk}

\begin{Def}[$ L^2- $ISS-Lyapunov function \cite{prieur2012iss}]
An $ L^2- $function, $\mathcal{L}, $ is said to be an ISS-Lyapunov function for the system \eqref{eq:LHSBLs} with the boundary conditions \eqref{eq:BCsLHSBLs} if there exist positive real constants $ \eta > 0 $, $ \xi > 0 $  and $ \beta > 0$ such that, for all continuous functions $ \Psi $, for all solutions of the system \eqref{eq:LHSBLs} satisfying the boundary conditions \eqref{eq:BCsLHSBLs}, and for all $ t \in [0, +\infty) $, 
\begin{equation}\label{AppxDtLyapunovfun}
\frac{d\mathcal{L}(W(\cdot,t))}{dt} \leq - \eta \mathcal{L}(W(\cdot,t)) + \frac{\beta}{\xi}\sup_{s \in [0,t]} \left(\int_{0}^{l}|\Psi(x,s)|^2dx\right).
\end{equation}
\end{Def}
At this point, inequalities that will be used in the proof of the main theorem are presented.
\begin{Prop}\label{Prop:Prop01}
	Let $ y , z \in \mathbb{R}^k$. Then,
	\begin{itemize}
		\item[a)] For any  matrix $ A \in \mathbb{R}^{k \times k} $, the following holds
		\begin{equation}\label{eq:QuadraticRelation1}
		- 2y^\top A(y - z) = - y^\top A y + z^\top A z - (y -z )^\top A (y -z ).
		\end{equation}
		\item[b)] For any positive semi-definite matrix $ B \in \mathbb{R}^{k \times k} $, there exist $ \xi > 0 $ such that 
		\begin{equation}\label{eq:QuadraticRelation2}
		\pm 2y^\top B z \leq \xi y^\top B y + \frac{1}{\xi}z^\top B z.
		\end{equation} 
	\end{itemize}
	\end{Prop}
\begin{Pf}
	a) For any quadratic form, 
	\begin{align*}
	(y -z )^\top A (y -z ) =&\; y^\top A y + z^\top A z - 2y^\top A z,\\
	=&\; - y^\top A y + z^\top A z - 2y^\top A z + 2y^\top A y,\\
	=&\; - y^\top A y + z^\top A z  + 2y^\top A(y - z).
	\end{align*}
	Thus equation \eqref{eq:QuadraticRelation1} is proved. 
	
	b) For any positive semi-definite quadratic form,  
	\begin{align*}
	0 \leq&\; \left(\sqrt{\xi}y \mp \frac{1}{\sqrt{\xi}}z\right)^\top B \left(\sqrt{\xi}y \mp \frac{1}{\sqrt{\xi}}z\right),\\
	=&\; {\xi} y^\top B y + \frac{1}{\xi}z^\top B z \mp 2y^\top B z.
	\end{align*}
	Thus, the inequality \eqref{eq:QuadraticRelation2} is obtained.
	\qed	
\end{Pf}
Further
\begin{Lem}\label{lem:Lemma01}
Consider the  $ L^2- $function, $ \mathcal{L}, $ defined by 
\begin{equation}\label{eq:LyapunovfunCh01}
\mathcal{L}(W(\cdot,t)) = \int_{0}^{l} W^{\top} P(x)W dx, \; t \in [0,+\infty),
\end{equation} which is an ISS-Lyapunov function for the system \eqref{eq:LHSBLs} with boundary conditions \eqref{eq:BCsLHSBLs}. 
Denote the smallest and largest eigenvalues of the diagonal matrix $ P(x) $ by $ \zeta $ and $ \beta $, respectively. Then, there exist a positive real constant $ \eta > 0 $, and for every $ W $, we have the inequality  
\begin{align}\label{eq:Lem01-cond02}
		\zeta \int_{0}^{l} |W|^2 dx &\leq \mathcal{L}(W(\cdot,t)) \leq \beta \int_{0}^{l} |W|^2 dx.
	\end{align}
\end{Lem}
The proof that $\mathcal{L}$ is indeed a Lyapunov function will be presented in Theorem \ref{thm:Ch01-01}. 
\begin{Rk}
In the paper \cite{bastin2017quadratic}, an explicit Lyapunov function is considered. The weight function, $P(x)$, used in the $ L^2- $ function \eqref{eq:LyapunovfunCh01} is in general form. Alternatively, the implicit weight function is defined by 
	\begin{equation*}
	P(x)= \text{diag}\{ {P^+}\exp(-\mu x) , {P^-}\exp(\mu x) \}, \quad \mu > 0, 
	\end{equation*}
	where $ P^+ \in \mathbb{R}_{+}^{m \times m} $ and $ P^- \in \mathbb{R}_{+}^{(k-m) \times (k-m)}$ are constant diagonal matrices.    
\end{Rk}
Now we are ready to present a proof of Lemma \ref{lem:Lemma01}:
\begin{Pf}
Since the diagonal matrix $ P(x)$ is a positive diagonal matrix, for all $x \in [0,l] $, and for all $ W \in \mathbb{R}^k $, we have  
\begin{equation}\label{eq:lllllll}
\zeta |W|^2 \leq W^{\top}P(x) W \leq \beta |W|^2.
\end{equation}
Consequently, the inequality \eqref{eq:Lem01-cond02} is obtained. 
\qed
\end{Pf}

\begin{Lem}[Gronwall's Lemma]\label{lem:Lemma02}
	Let $ y \in C^1([0, +\infty)) $, $ z \in \mathbb{R} $, $ a \in \mathbb{R}^+ $, and 
	\begin{equation*}
		y'(t) \leq - a y(t) + z,\; y(0) = c \geq 0, \; t \geq 0.
	\end{equation*}
	Then 
	\begin{equation*}
		y(t) \leq \left(c - \frac{z}{a}\right) e^{-at} + \frac{z}{a}, \; t \geq 0.
	\end{equation*}
\end{Lem}
\begin{Pf}
	For the proof, see Lemma 1.1.1 in \cite{lakshmikantham1989stability} by considering constants $ a $ and $ z $. \qed
\end{Pf}

\begin{Thm}[Stability]\label{thm:Ch01-01}
Assume the system \eqref{eq:LHSBLs} with the boundary conditions \eqref{eq:BCsLHSBLs} satisfies Assumption \ref{Asmpt:LHSBLs-Asmpt01}. Let $ \xi $ be any positive real number. Define a weight function by $P(x)= \text{diag}\{ {P^+}(x) , {P^-}(x) \}$, where $ P^+(x) \in \mathbb{R}_{+}^{m \times m}$ and $P^-(x) \in \mathbb{R}_{+}^{(k-m) \times (k-m)}$. Assume that the matrix 
\begin{equation}\label{eq:Theorem01-cond01}
- \Lambda(x) P'(x) - \Lambda'(x) P(x) + \Pi^{\top}(x)P(x) + P(x) \Pi(x) - \xi P(x),
\end{equation}
is positive definite for all $ x \in [0,l]$  and the matrix  
\begin{equation}\label{eq:Theorem01-cond02}
\begin{bmatrix} \Lambda^+(l)P^+(l) & 0 \\ 0 &  \Lambda^-(0) P^-(0) \end{bmatrix} - K^{\top} \begin{bmatrix} \Lambda^+(0)P^+(0)  & 0 \\ 0 & - \Lambda^-(l) P^-(l) \end{bmatrix} K,
\end{equation}
is positive semi-definite. Then the  $ L^2- $function, $ \mathcal{L}, $ in Equation \eqref{eq:LyapunovfunCh01} 
is an ISS-Lyapunov function for the system \eqref{eq:LHSBLs} with boundary conditions \eqref{eq:BCsLHSBLs}. Moreover, the steady-state $ W(x,t) \equiv 0 $ of the system \eqref{eq:LHSBLs} with boundary conditions \eqref{eq:BCsLHSBLs} is ISS in $L^2-$norm with respect to the disturbance function $ \Psi $.
\end{Thm}

\begin{Pf}
It suffices to show that the $ L^2-$function defined by Equation \eqref{eq:LyapunovfunCh01} is an ISS-Lyapunov function. Thus, the time derivative of the candidate ISS-Lyapunov function \eqref{eq:LyapunovfunCh01} is computed by using the system \eqref{eq:LHSBLs}, the boundary conditions \eqref{eq:BCsLHSBLs}, and the initial boundary compatibility conditions \eqref{eq:CCs01LHSBLs} as follows:
\begin{align*}
	\frac{d\mathcal{L}(W(\cdot,t))}{dt} 
	= &\; \int_{0}^{l} \left(\partial_t{W^{\top}}P(x)W + {W^{\top}} P(x)\partial_t W\right) dx,\\
	= &\; \int_{0}^{l}\left( \left(-\Lambda (x) \partial_x W  -\Pi(x) W + \Psi(x,t)\right)^{\top}P(x)W \right.\\
	&\; \left.+ {W^{\top}} P(x)\left(-\Lambda (x) \partial_x W  -\Pi(x) W + \Psi(x,t)\right) \right) dx.
\end{align*}
Rewriting the above equation, one obtains:
\begin{eqnarray*}
\frac{d\mathcal{L}(W(\cdot,t))}{dt} & = &\; -\int_{0}^{l} \left(\partial_x{W^{\top}}\Lambda(x) P(x)  W + {W^{\top}}\Lambda(x) P(x) \partial_x W \right) \; dx \\
&&\; - \int_{0}^{l} {W^{\top}} \left(\Pi^{\top}(x)P(x) + P(x) \Pi(x)\right)Wdx + 2 \int_{0}^{l} {W^{\top}} P(x) \Psi(x,t) \; dx,\\
&&\; - \int_{0}^{l} \partial_x \left( {W^{\top}}\Lambda(x) P(x) W \right)dx  + 2\int_{0}^{l} {W^{\top}} P(x) \Psi(x,t)\; dx\\
&&\;- \int_{0}^{l} {W^{\top}} \left( - \Lambda(x) P'(x) - \Lambda'(x) P(x) + \Pi^{\top}(x)P(x) + P(x)\Pi(x) \right)W\;dx,\\
	&= &\; - \left[ {W^{\top}}\Lambda(x) P(x) W \right]_0^l + 2\int_{0}^{l} {W^{\top}} P(x) \Psi(x,t) \; dx  \\
	&& - \int_{0}^{l} {W^{\top}} \left(- \Lambda(x) P^\prime(x) - \Lambda^\prime(x) P(x) + \Pi^{\top}(x)P(x) + P(x)\Pi(x) \right)W\; dx.
\end{eqnarray*}
The first term in the above equation is treated by substituting boundary conditions \eqref{eq:BCsLHSBLs} and by assumption in Theorem \ref{thm:Ch01-01} for the matrix in  Equation \eqref{eq:Theorem01-cond02} to obtain:
\begin{align*}
	- \left[ {W^{\top}}\Lambda(x) P(x) W \right]_0^l 
	= &\; -\begin{bmatrix} W^+(l,t) \\ W^-(l,t) \end{bmatrix}^{\top} \begin{bmatrix} \Lambda^+(l)P^+(l) & 0 \\ 0 & - \Lambda^-(l) P^-(l) \end{bmatrix} \begin{bmatrix} W^+(l,t) \\ W^-(l,t) \end{bmatrix}\\
	&\; +\begin{bmatrix} W^+(0,t) \\ W^-(0,t) \end{bmatrix}^{\top} \begin{bmatrix} \Lambda^+(0)P^+(0)  & 0 \\ 0 & - \Lambda^-(0) P^-(0) \end{bmatrix} \begin{bmatrix} W^+(0,t) \\ W^-(0,t) \end{bmatrix},\\
	= &\; -\begin{bmatrix} W^+(l,t) \\ W^-(0,t) \end{bmatrix}^{\top} \begin{bmatrix} \Lambda^+(l)P^+(l) & 0 \\ 0 &  \Lambda^-(0) P^-(0) \end{bmatrix} \begin{bmatrix} W^+(l,t) \\ W^-(0,t) \end{bmatrix}\\
	&\; +\begin{bmatrix} W^+(0,t) \\ W^-(l,t) \end{bmatrix}^{\top} \begin{bmatrix} \Lambda^+(0)P^+(0)  & 0 \\ 0 & - \Lambda^-(l) P^-(l) \end{bmatrix} \begin{bmatrix} W^+(0,t) \\ W^-(l,t) \end{bmatrix},\\
	= &\; -\begin{bmatrix} W^+(l,t) \\ W^-(0,t) \end{bmatrix}^{\top} \begin{bmatrix} \Lambda^+(l)P^+(l) & 0 \\ 0 &  \Lambda^-(0) P^-(0) \end{bmatrix} \begin{bmatrix} W^+(l,t) \\ W^-(0,t) \end{bmatrix}\\
	&\; +\begin{bmatrix} W^+(l,t) \\ W^-(0,t) \end{bmatrix}^{\top} K^{\top} \begin{bmatrix} \Lambda^+(0)P^+(0)  & 0 \\ 0 & - \Lambda^-(l) P^-(l) \end{bmatrix} K \\
	&\; \begin{bmatrix} W^+(l,t) \\ W^-(0,t)\end{bmatrix} \le 0.
\end{align*}

Thus, 
\begin{align}
	\frac{d\mathcal{L}(W(\cdot,t))}{dt} = &\;-\int_{0}^{l} {W^{\top}} \left( - \Lambda(x) P'(x) - \Lambda'(x) P(x) + \Pi^{\top}(x)P(x) + P(x)\Pi(x) \right)W dx  \nonumber\\
	&\; + 2\int_{0}^{l} {W^{\top}} P(x) \Psi(x,t) dx. \label{eq:lll}
\end{align}

By Proposition \ref{Prop:Prop01}, and using the inequality \eqref{eq:Lem01-cond02}, the time derivative of the candidate ISS-Lyapunov function \eqref{eq:lll} is estimated as 
\begin{eqnarray*}
		\frac{d\mathcal{L}(W(\cdot,t))}{dt} &\leq &\; -\int_{0}^{l} {W^{\top}} \left[-\Lambda(x) P'(x) - \Lambda'(x) P(x) +\Pi^{\top}(x)P(x) + P(x)\Pi(x) \right]Wdx  \\
		&&\;  + \xi\int_{0}^{l} {W^{\top}} P(x) W dx + \frac{1}{\xi}\int_{0}^{l} {\Psi^{\top}(x,t)} P(x) \Psi(x,t) dx,\\		
		&\leq &\; -\int_{0}^{l} {W^{\top}} \left[-\Lambda(x) P'(x) - \Lambda'(x) P(x) +\Pi^{\top}(x)P(x) + P(x)\Pi(x)  - \xi P(x) \right]Wdx  \nonumber\\
		&&\;  + \frac{\beta}{\xi}\int_{0}^{l} |\Psi(x,t)|^2 dx.	
\end{eqnarray*}
Hence
\begin{equation}\label{eq:llll}
\frac{d\mathcal{L}(W(\cdot,t))}{dt} \le - \int_{0}^{l} {W^{\top}} Q(x)Wdx + \frac{\beta}{\xi}\sup_{s \in [0,t]} \left(\int_{0}^{l}|\Psi(x,s)|^2dx\right).
\end{equation}
	where $ Q(x) = - \Lambda(x) P'(x) - \Lambda'(x) P(x) +\Pi^{\top}(x)P(x) + P(x)\Pi(x) - \xi P(x) $.
	
	Also using assumptions in Theorem \ref{thm:Ch01-01} for the matrix in Equation \eqref{eq:Theorem01-cond01}, $ Q $, there exists $ \eta > 0 $ such that
	$ W^{\top}Q(x) W  \geq \eta W^{\top}P(x) W $ for all $ x \in [0,l]$, we obtain the inequality \eqref{eq:Lem01-cond01} below:
\begin{align}
		\frac{d\mathcal{L}(W(\cdot,t))}{dt} &\leq - \eta \mathcal{L}(W(\cdot,t)) + \frac{\beta}{\xi}\sup_{s \in [0,t]} \left(\int_{0}^{l}|\Psi(x,s)|^2dx\right).\label{eq:Lem01-cond01}
	\end{align}
We now have 
\begin{align}
	\mathcal{L}(W(\cdot,t)) \leq&\; e^{-\eta t} \left(\mathcal{L}(W(\cdot,0)) - \frac{\beta}{\eta \xi}\sup_{s \in [0,t]} \left(\int_{0}^{l}|\Psi(x,s)|^2dx\right) \right)\nonumber \\
	&\; + \frac{\beta}{\eta \xi}\sup_{s \in [0,t]} \left(\int_{0}^{l}|\Psi(x,s)|^2dx\right),\nonumber\\
	&\;\leq e^{-\eta t} \mathcal{L}(W(\cdot,0)) + \frac{\beta}{\eta \xi}\sup_{s \in [0,t]} \left(\int_{0}^{l}|\Psi(x,s)|^2dx\right),\;t \geq 0.\label{eq:llllll}
\end{align}
Now, we apply inequality \eqref{eq:Lem01-cond02} in inequality \eqref{eq:llllll} to obtain
\begin{align}
	\zeta{\|W(\cdot, t) \|}_{L^2((0,l);\mathbb{R}^k)}^2 \leq \beta {e}^{-\eta t}{\|W_0 \|}_{L^2((0,l);\mathbb{R}^k)}^2 + \frac{\beta}{\eta \xi}\sup_{s \in [0,t]} \left(\int_{0}^{l}|\Psi(x,s)|^2dx\right),\;t \geq 0.
\end{align}
Let $ C = \beta/\zeta $. Then the condition for exponential stability Equation \eqref{eq:ExponStabCondChap1} is satisfied. This concludes the proof of Theorem \ref{thm:Ch01-01}.
\qed	 
\end{Pf}

Now consider a $ k\times k $ uniform linear hyperbolic system of balance laws with additive disturbances: 
\begin{equation}\label{eq:LHSBLs2}
\partial_t W + \Lambda \partial_x W + \Pi W  = \Psi(x,t),
\end{equation}
where $\Lambda = \text{diag}\{{\Lambda^+}, -{\Lambda^-} \} $, with $ {\Lambda^+} \in \mathbb{R}_{+}^{m \times m}$, and $ {\Lambda^-} \in \mathbb{R}_{+}^{(k-m) \times (k-m)}$, is a diagonal matrix, $ \Pi $ is a constant real matrix in $ \mathbb{R}^{k\times k} $, and $ \Psi $ is a vector in $ \mathbb{R}^k $.
\begin{Asmpt}\label{Asmpt:LHSBLs-Asmpt02}
For all $ x \in [0,l]$, and $ t \in [0,+\infty)$, assume that assumption (ii) for $ \Psi $  and (iii) stated in Assumption \ref{Asmpt:LHSBLs-Asmpt01} still hold.
\end{Asmpt}

\begin{Cc}\label{col:Col01}
Assume the system \eqref{eq:LHSBLs2} with the boundary conditions \eqref{eq:BCsLHSBLs} satisfies Assumption \ref{Asmpt:LHSBLs-Asmpt02}. Let $ \xi $ be any positive real number. Assume that the matrix 
\begin{equation}\label{eq:Col01-cond01}
- \Lambda P'(x) + \Pi^{\top}(x)P(x) + P(x) \Pi(x) - \xi P(x),
\end{equation}
is positive definite for all $ x \in [0,l]$, and the matrix  
\begin{equation}\label{eq:Col01-cond02}
\begin{bmatrix} \Lambda^+P^+(l) & 0 \\ 0 &  \Lambda^- P^-(0) \end{bmatrix} - K^{\top} \begin{bmatrix} \Lambda^+P^+(0)  & 0 \\ 0 & - \Lambda^- P^-(l) \end{bmatrix} K,
\end{equation}
is positive semi-definite. Then the $ L^2- $function, $ \mathcal{L} $, defined by Equation \eqref{eq:LyapunovfunCh01} is an ISS-Lyapunov function for the system \eqref{eq:LHSBLs2} with boundary conditions \eqref{eq:BCsLHSBLs}. Moreover, the steady-state $ W(x,t) \equiv 0 $ of the system \eqref{eq:LHSBLs2} with boundary conditions \eqref{eq:BCsLHSBLs} is ISS in $L^2-$norm with respect to disturbance function $ \Psi $.
\end{Cc}

\section{Numerical boundary feedback and ISS}\label{sec:sec02}
In this section a numerical approach to analyse the non-uniform balance laws is presented. For non-uniform linear hyperbolic systems of balance laws with additive disturbances in one spatial dimension the finite volume method is applied (see \cite{leveque2002finite}). 

Consider a uniform grid and denote grid points along the $ x- $ and $ t- $directions by  
\[ x_{j - \frac{1}{2}} = j {\Delta x},\; j = 0, \dots, J,\quad  t^n = n {\Delta t},\; n = 0, \dots, N, \]
respectively, where $ {\Delta x} = l/J $ and $ {\Delta t} = T/N $ denote step sizes, and $ x_{- \frac{1}{2}} = 0 $ and $ x_{J - \frac{1}{2}} = l $ denote the left and right boundary points, respectively. Let $ x_{j} = \left(j + \frac{1}{2}\right) {\Delta x}, \; j = 0, \dots, J-1 $ denote cell centres. Approximate the $ j-$th cell average at time $ t^n$ of the state variables $ W $ over each grid cell $ \left(x_{j - \frac{1}{2}}, x_{j + \frac{1}{2}} \right) $ by:
\begin{equation}\label{eq:CellAvg}
W_{j}^{n} = \frac{1}{\Delta x}\int_{x_{j - \frac{1}{2}}}^{x_{j + \frac{1}{2}}} W(x,t^n) dx, \; j = 0, \dots, J-1.
\end{equation}
With the numerical approximation $  W_{j}^{n} $, and for $ {\Delta t} \rightarrow 0$, the following operator-splitting technique is applied to the system \eqref{eq:LHSBLs}:
\begin{subequations}\label{eq:LHSBLs-Splitted}
\begin{align}
\partial_t W + \Lambda (x) \partial_x W =&\;\Psi(x,t),\label{eq:LHSBLs-Splitted-01}\\
\partial_t W + \Pi(x) W =&\; 0,\; x \in [0,l],\; t \in [0,+\infty).\label{eq:LHSBLs-Splitted-02}
\end{align}
\end{subequations}
Thus the system \eqref{eq:LHSBLs} is discretised by applying an explicit Euler scheme for temporal derivatives, an upwind scheme for spatial derivatives and centred discretisation for coefficients, source terms and additive disturbances as: for $  n = 0, \dots , N-1 $, $ j = 0, \dots, J-1 $, 
\begin{subequations}\label{eq:DiscLHSBLs}
\begin{align}
	&\begin{bmatrix} \widetilde{W^+}_{j}^{n}\\\widetilde{W^-}_{j}^{n} \end{bmatrix} =  \begin{bmatrix} {W^+}_{j}^{n}\\{W^-}_{j}^{n} \end{bmatrix}  - {\frac{\Delta t}{\Delta x}} \begin{bmatrix} {\Lambda^+}_{j-1} & 0 \\ 0 & - {\Lambda^-}_{j+1} \end{bmatrix} \begin{bmatrix}  {W^+}_{j}^n - {W^+}_{j-1}^n \\ {W^-}_{j+1}^n - {W^-}_{j}^n \end{bmatrix}+ {\Delta t} \begin{bmatrix} {\Psi^+}_j^n \\ {\Psi^-}_j^n \end{bmatrix}, \label{eq:DiscLHSBLs01}\\
	&\begin{bmatrix} {W^+}_{j}^{n+1}\\{W^-}_{j}^{n+1} \end{bmatrix} =  
	\begin{bmatrix} \widetilde{W^+}_{j}^{n}\\\widetilde{W^-}_{j}^{n} \end{bmatrix} - {\Delta t}\Pi_j \begin{bmatrix} \widetilde{W^+}_{j}^{n}\\\widetilde{W^-}_{j}^{n}  \end{bmatrix}.\label{eq:DiscLHSBLs02}
\end{align}
\end{subequations}  
The initial condition \eqref{eq:ICLHSBLs} is discretised as 
\begin{equation}\label{eq:DiscICLHSBLs}
W_j^0 = W_{0,j}, \quad j = 0, \dots,J-1,
\end{equation}
and the discretisation of the boundary conditions \eqref{eq:BCsLHSBLs} is
\begin{equation}\label{eq:DiscBCsLHSBLs}
\begin{bmatrix} {W^+}_{-1}^{n+1}\\{W^-}_{J}^{n+1} \end{bmatrix} =  K \begin{bmatrix} {W^+}_{J-1}^{n+1}\\{W^-}_{0}^{n+1} \end{bmatrix},\quad n = 0, \dots , N-1.
\end{equation}
Furthermore, the discretisation of the zero-order initial boundary compatibility conditions \eqref{eq:CCs01LHSBLs} can be written as
\begin{equation}\label{eq:DiscCCs01LHSBLs}
\begin{bmatrix} {W^+}_{-1}^{0}\\{W^-}_{J}^{0} \end{bmatrix} =  K \begin{bmatrix} {W^+}_{J-1}^{0}\\{W^-}_{0}^{0} \end{bmatrix}. 
\end{equation}

Our aim is to investigate conditions for numerical boundary feedback stabilisation. For this reason, the definition of discrete ISS follows: 
\begin{Def}[Discrete ISS]
The steady-state $ W_j^n \equiv 0,\; j = 0, \dots, J-1,\; n = 0, \dots, N-1 $ of the discretised system \eqref{eq:DiscLHSBLs}, with boundary conditions \eqref{eq:DiscBCsLHSBLs} is discrete ISS in $L^2-$norm with respect to discrete disturbance function $ \Psi_j^n$ if there exist positive real constants $\eta > 0$, $ \xi > 0 $ and $C>0$ such that, for every initial condition $W_j^0 \in L^2((x_{j-\frac{1}{2}},x_{j+\frac{1}{2}});\mathbb{R}^k)$ satisfying the compatibility condition \eqref{eq:DiscCCs01LHSBLs}, the  $L^2-$solution to the discretised system \eqref{eq:DiscLHSBLs} with initial conditions \eqref{eq:DiscICLHSBLs} and boundary conditions \eqref{eq:DiscBCsLHSBLs} satisfies
\begin{equation}\label{eq:DiscExponStabCondChap1}
	{\Delta x}{\sum_{j=0}^{J-1}} |W_j^{n+1}|^2 \leq Ce^{-\eta t^{n+1}}{\Delta x}{\sum_{j=0}^{J-1}} |W_j^0|^2 +  \frac{C}{\eta}\left(\frac{1}{\xi} + {\Delta t}\right) \sup_{0 \leq s \leq n}\left({\Delta x}\sum_{j=0}^{J-1} |\Psi_j^s|^2\right). 
	\end{equation}
\end{Def}

\begin{Def}[A discrete $ L^2 -$ISS-Lyapunov function]\label{def:DiscISSLyafun}
A discrete $ L^2 $ function, $ \mathcal{L}^n,\; n = 0, \dots, N-1$ is said to be a discrete ISS-Lyapunov function for the discretised system \eqref{eq:DiscLHSBLs} with boundary conditions \eqref{eq:DiscBCsLHSBLs} if there exist positive real constants $ \eta > 0 $, $ \xi > 0 $  and $ \beta > 0 $ such that, for all discrete functions $ \Psi_j^n,\;j=0,\dots, J-1 $, for all solutions of the system \eqref{eq:DiscLHSBLs} satisfying boundary conditions \eqref{eq:DiscBCsLHSBLs},
\begin{equation}\label{AppxDtDiscLyapunovfun}
	\frac{\mathcal{L}^{n+1} - \mathcal{L}^{n}}{\Delta t} \leq - \eta \mathcal{L}^{n} + \beta \left(\frac{1}{\xi} + {\Delta t}\right)  \sup_{0 \leq s \leq n}\left({\Delta x}\sum_{j=0}^{J-1} |\Psi_j^s|^2\right).
\end{equation}
\end{Def}
Before we proceed with the main result of this section, we present the following preliminary results.
\begin{Lem}
Let the discrete $ L^2-$function defined by 
\begin{equation}\label{eq:DiscLyapunovfunCh01}
	\mathcal{L}^n = {\Delta x }\sum_{j=0}^{J-1} {W_j^n}^{\top} P_j W_j^n, \quad n = 0, \dots , N-1,
\end{equation} 
be a discrete ISS-Lyapunov function for the system \eqref{eq:DiscLHSBLs} with boundary conditions \eqref{eq:DiscBCsLHSBLs}. Define a discrete weight function by $P_j= \text{diag}\{ P_j^+ , P_j^-\}$, where $ P_j^+$ and $P_j^-$ denote the first $ m $ and the last $ k-m $  diagonal entries, respectively, for $j = 0, \dots, J-1$ with the smallest and largest eigenvalue of $ P_j,\ j = 0, \dots, J-1 $ denoted by $\displaystyle \zeta  = \min_{0 \leq j \leq J-1}P_j $ and $\displaystyle \beta  = \max_{0 \leq j \leq J-1}P_j $, respectively.  Then the following inequality holds  
\begin{align}
	\zeta {\Delta x }\sum_{j=0}^{J-1} |W_j^n|^2  &\leq \mathcal{L}^{n} \leq \beta  {\Delta x }\sum_{j=0}^{J-1} |W_j^n|^2.\label{eq:DiscLem01-cond02}
\end{align}
\end{Lem}
\begin{Pf}
Consider the positive diagonal matrix $ P_j $, then we have  
\begin{equation}\label{eq:PostiveSD02}
	\zeta |W_j^n|^2 \leq {W}_{j}^{n \top} P_j {W}_{j}^{n} \leq \beta |W_j^n|^2,\; j = 0, \dots, J-1,\; n = 0, \dots, N-1.
\end{equation}
Thus, the inequality \eqref{eq:DiscLem01-cond02} can be obtained from the inequality \eqref{eq:PostiveSD02}. 
\end{Pf}

\begin{Lem}\label{lem:DiscLemma2}
Let $ a > 0 $ and $ z \in \mathbb{R} $. Suppose for discrete functions $ y^{n}, \; n = 0 , \dots, N-1  $, 
\begin{equation}\label{eq:DiscDInequality}
	\frac{y^{n+1} - y^{n}}{\Delta t} \leq - a y^{n} + z,\; y^{0} = c. 
\end{equation}
Then
\begin{equation}\label{eq:DiscInequality}
	y^{n+1} \leq \left(c - \frac{z}{a}\right) \left(1 - a{\Delta t} \right)^{n+1} + \frac{z}{a}, \; n = 0 , \dots, N-1.
\end{equation}
\end{Lem}

\begin{Pf}
	The inequality \eqref{eq:DiscDInequality}, by applying recursion, can be expressed as
	\begin{equation}\label{eq:DiscDInequality1}
		y^{n+1} \leq c \left(1 - a{\Delta t} \right)^{n+1} + z{\Delta t} \sum_{r = 0}^{n} \left(1 - a{\Delta t} \right)^{r}, \; \; n = 0 , \dots, N-1.
	\end{equation}
	For sufficiently small $ {\Delta t}$, $ 0 < 1 - a{\Delta t} < 1 $, then the inequality \eqref{eq:DiscDInequality1} implies the inequality \eqref{eq:DiscInequality}.
	\qed
\end{Pf}

\begin{Thm}[Stability]\label{thm:ThmChap01-02}
Assume the system \eqref{eq:DiscLHSBLs} with boundary conditions \eqref{eq:DiscBCsLHSBLs} satisfies Assumption \eqref{Asmpt:LHSBLs-Asmpt01} for system \eqref{eq:DiscLHSBLs}. Let $ T > 0 $ be fixed and the CFL condition, $ \frac{\Delta t}{\Delta x} \max_{\stackrel{1\leq i \leq k}{0 \leq j \leq J-1}} |\lambda_{i,j}| \leq 1 $ hold. Let $ \xi $ be any positive real number. Define a discrete weight function by $P_j= \text{diag}\{ P_j^+ , P_j^-\}$, where $ P_j^+$ and $P_j^-$ denote the first $ m $ and the last $ k-m $  diagonal entries, respectively, for $j = 0, \dots, J-1$. Assume that the matrix  
\begin{multline}\label{eq:DiscTheorem01-cond00}
-\left(1 + \xi {\Delta t} \right)\begin{bmatrix} {\Lambda_{j-1}^+}\left(\frac{P_{j+1}^+ - P_{j}^+}{\Delta x} \right)  & 0 \\ 0 &  -{\Lambda_{j+1}^-}\left(\frac{P_{j}^- - P_{j-1}^-}{\Delta x} \right)  \end{bmatrix}\\
-\left(1 + \xi {\Delta t} \right) \begin{bmatrix} \left(\frac{{\Lambda_{j}^+} - {\Lambda_{j-1}^+}}{\Delta x}\right)P_{j+1}^+ & 0 \\ 0 &   -\left(\frac{{\Lambda_{j+1}^-} - {\Lambda_{j}^-}}{\Delta x}\right)P_{j-1}^- \end{bmatrix} -\xi P_j,
\end{multline}
is positive definite for all $ j = 0, \dots, J-1$ and the matrices
\begin{equation}\label{eq:DiscTheorem01-cond01}
	P_j\Pi_j + {\Pi_j}^{\top} P_j - {\Delta t}{\Pi_j}^{\top}P_j{\Pi_j},
\end{equation}
and 
\begin{equation}\label{eq:DiscTheorem01-cond02}
	\begin{bmatrix} {\Lambda^+}_{J-1}P_{J}^+ & 0 \\ 0 &  {\Lambda^-}_{0} P_{-1}^- \end{bmatrix} - K^{\top}\begin{bmatrix}  {\Lambda^+}_{-1}P_{0}^+ & 0 \\ 0 &  {\Lambda^-}_{J}P_{J-1}^- \end{bmatrix}K,
\end{equation}
are positive semi-definite for all  $ j = 0, \dots, J-1$. Then the discrete $ L^2-$function defined by Equation \eqref{eq:DiscLyapunovfunCh01}
is a discrete ISS-Lyapunov function for system \eqref{eq:DiscLHSBLs} with boundary conditions \eqref{eq:DiscBCsLHSBLs}. Moreover, the steady-state $ W_j^n \equiv 0,\; j = 0, \dots, J-1,\; n = 0, \dots, N-1 $ of system \eqref{eq:DiscLHSBLs} with boundary conditions \eqref{eq:DiscBCsLHSBLs} is discrete ISS in $L^2-$norm with respect to discrete disturbance function $ \Psi_j^n,\;j=0,\dots, J-1 $.
\end{Thm}
\begin{Pf}
We begin the proof by approximating the time derivative of the candidate ISS-Lyapunov function defined by Equation \eqref{eq:LyapunovfunCh01}. It can be expressed as 
\begin{equation}\label{eq:DtDiscLyapunovfun}
	\frac{\mathcal{L}^{n+1} - \mathcal{L}^{n}}{\Delta t} = \frac{\mathcal{L}^{n+1} - \widetilde{\mathcal{L}}^{n}}{\Delta t} + \frac{\widetilde{\mathcal{L}}^{n} - \mathcal{L}^{n}}{\Delta t},
\end{equation}
where \[\widetilde{\mathcal{L}}^{n} = {\Delta x }\sum_{j=0}^{J-1} \widetilde{W}_j^{n\top} P_j \widetilde{W}_j^n, \quad n = 0, \dots , N-1.\]

Consider the first term on the RHS of equation \eqref{eq:DtDiscLyapunovfun} and then by system \eqref{eq:DiscLHSBLs02}, we have 
\begin{align}
\frac{\mathcal{L}^{n+1} - \widetilde{\mathcal{L}}^{n}}{\Delta t} =&\; \frac{\Delta x }{\Delta t}\sum_{j=0}^{J-1}\left(W_j^{n+1 \top} P_j W_j^{n+1} - \widetilde{W}_j^{n \top} P_j \widetilde{W}_j^n\right),\nonumber \\
	= &\; \frac{\Delta x}{\Delta t} \sum_{j=0}^{J-1} \begin{aligned}[t]
	&\left( \widetilde{W}_j^{n \top} P_j \widetilde{W}_j^n - {\Delta t} \widetilde{W}_j^{n \top} P_j\Pi_j \widetilde{W}_j^n  - {\Delta t} \widetilde{W}_j^{n \top}{\Pi_j}^{\top} P_j \widetilde{W}_j^n \right.\\
	&\left.  + \left({\Delta t}\right)^2 \widetilde{W}_j^{n \top} {\Pi_j}^{\top} P_j {\Pi_j} \widetilde{W}_j^n  - \widetilde{W}_j^{n \top} P_j \widetilde{W}_j^n\right),\end{aligned}\nonumber\\
	=&\; - {\Delta x} \sum_{j=0}^{J-1} \widetilde{W}_j^{n \top} \left(P_j\Pi_j + {\Pi_j}^{\top} P_j - {\Delta t}{\Pi_j}^{\top} P_j {\Pi_j} \right) \widetilde{W}_j^n \leq 0, \label{eq:DtDiscLyapunovfun01}	
\end{align} 
where $n = 0, \dots, N-1$ and we used the assumption in Theorem \ref{thm:ThmChap01-02} in the last step. 

We now analyse the second term on the RHS of the equation \eqref{eq:DtDiscLyapunovfun} as
\begin{equation}
\frac{\widetilde{\mathcal{L}}^n-\mathcal{L}^{n}}{\Delta t} =\; \frac{\Delta x}{\Delta t} \sum_{j=0}^{J-1} \left(\widetilde{W}_j^{n \top} P_j \widetilde{W}_j^{n} - W_j^{n \top} P_j W_j^n \right),\; n = 0, \dots, N-1, \label{eq:DtDiscLyapunovfun02}
\end{equation}
where 
\begin{eqnarray}
\widetilde{W}_j^{n \top} P_j \widetilde{W}_j^{n} &=& W_j^{n \top} P_j W_j^n - 2\begin{bmatrix} {W^+}_{j}^{n}\\{W^-}_{j}^{n} \end{bmatrix}^{\top}  \begin{bmatrix} {\frac{\Delta t}{\Delta x}}{\Lambda_{j-1}^+}P_j^+ & 0 \\ 0 & - {\frac{\Delta t}{\Delta x}}{\Lambda_{j+1}^-}P_j^- \end{bmatrix} \begin{bmatrix}  {W^+}_{j}^n - {W^+}_{j-1}^n \\ {W^-}_{j+1}^n - {W^-}_{j}^n \end{bmatrix}\nonumber\\
&& +\begin{bmatrix}  {W^+}_{j}^n - {W^+}_{j-1}^n \\ {W^-}_{j+1}^n - {W^-}_{j}^n \end{bmatrix}^{\top} \begin{bmatrix} \left({\frac{\Delta t}{\Delta x}}{\Lambda_{j-1}^+}\right)^2P_j^+ & 0 \\ 0 & \left(- {\frac{\Delta t}{\Delta x}}{\Lambda_{j+1}^-}\right)^2 P_j^- \end{bmatrix}   \begin{bmatrix}  {W^+}_{j}^n - {W^+}_{j-1}^n \\ {W^-}_{j+1}^n - {W^-}_{j}^n \end{bmatrix} \nonumber\\
&& + 2 {\Delta t} \begin{bmatrix} {W^+}_{j}^{n}\\{W^-}_{j}^{n} \end{bmatrix}^\top \begin{bmatrix} \left(I_{m} - {\frac{\Delta t}{\Delta x}}{\Lambda^+}_{j-1}\right) P_j^+ & 0 \\ 0 & \left(I_{k-m} - {\frac{\Delta t}{\Delta x}}{\Lambda^-}_{j+1}\right) P_j^- \end{bmatrix}\begin{bmatrix} {\Psi^+}_j^n \\ {\Psi^-}_j^n \end{bmatrix}\label{eq:DtDiscLyapunovfun02-001}  \\
&& + 2{\Delta t} \begin{bmatrix} {W^+}_{j-1}^n \\ {W^-}_{j+1}^n \end{bmatrix}^{\top}  \begin{bmatrix} {\frac{\Delta t}{\Delta x}}{\Lambda_{j-1}^+}P_j^+ & 0 \\ 0 & {\frac{\Delta t}{\Delta x}}{\Lambda_{j+1}^-}P_j^- \end{bmatrix}\begin{bmatrix} {\Psi^+}_j^n \\ {\Psi^-}_j^n \end{bmatrix}  + \left({\Delta t}\right)^2 {\Psi_j^n}^{\top}P_j\Psi_j^n. \nonumber
\end{eqnarray}

By using Proposition \ref{Prop:Prop01} and the CFL condition in equation \eqref{eq:DtDiscLyapunovfun02-001}, for all $n = 0, \dots, N-1 $, $j = 0, \dots, J-1 $, we obtain  
\begin{eqnarray*}
\widetilde{W}_j^{n \top} P_j \widetilde{W}_j^{n} &\leq & W_j^{n \top} P_j W_j^n - 
\begin{bmatrix} {W^+}_{j}^{n}\\{W^-}_{j}^{n} \end{bmatrix}^{\top}  \begin{bmatrix} {\frac{\Delta t}{\Delta x}}{\Lambda_{j-1}^+}P_j^+ & 0 \\ 0 &  {\frac{\Delta t}{\Delta x}}{\Lambda_{j+1}^-}P_j^- \end{bmatrix} \begin{bmatrix} {W^+}_{j}^{n}\\{W^-}_{j}^{n} \end{bmatrix} \\
&&+ \begin{bmatrix} {W^+}_{j-1}^{n}\\{W^-}_{j+1}^{n} \end{bmatrix}^{\top}  \begin{bmatrix} {\frac{\Delta t}{\Delta x}}{\Lambda_{j-1}^+}P_j^+ & 0 \\ 0 &  {\frac{\Delta t}{\Delta x}}{\Lambda_{j+1}^-}P_j^- \end{bmatrix} \begin{bmatrix} {W^+}_{j-1}^{n}\\{W^-}_{j+1}^{n} \end{bmatrix}\\
&& -\begin{bmatrix}  {W^+}_{j}^n - {W^+}_{j-1}^n \\ {W^-}_{j+1}^n - {W^-}_{j}^n \end{bmatrix}^{\top} \begin{bmatrix} {\frac{\Delta t}{\Delta x}}{\Lambda_{j-1}^+}P_j^+ & 0 \\ 0 &  {\frac{\Delta t}{\Delta x}}{\Lambda_{j+1}^-}P_j^- \end{bmatrix}   \begin{bmatrix}  {W^+}_{j}^n - {W^+}_{j-1}^n \\ {W^-}_{j+1}^n - {W^-}_{j}^n \end{bmatrix}\\
&& +\begin{bmatrix}  {W^+}_{j}^n - {W^+}_{j-1}^n \\ {W^-}_{j+1}^n - {W^-}_{j}^n \end{bmatrix}^{\top} \begin{bmatrix} \left({\frac{\Delta t}{\Delta x}}{\Lambda_{j-1}^+}\right)^2P_j^+ & 0 \\ 0 & \left(- {\frac{\Delta t}{\Delta x}}{\Lambda_{j+1}^-}\right)^2 P_j^- \end{bmatrix}   \begin{bmatrix}  {W^+}_{j}^n - {W^+}_{j-1}^n \\ {W^-}_{j+1}^n - {W^-}_{j}^n \end{bmatrix}\\
&& +  {\Delta t} \xi \begin{bmatrix} {W^+}_{j}^{n}\\{W^-}_{j}^{n} \end{bmatrix}^\top \begin{bmatrix} \left(I_{m} - {\frac{\Delta t}{\Delta x}}{\Lambda^+}_{j-1}\right) P_j^+ & 0 \\ 0 & \left(I_{k-m} - {\frac{\Delta t}{\Delta x}}{\Lambda^-}_{j+1}\right) P_j^- \end{bmatrix}\begin{bmatrix} {W^+}_{j}^{n}\\{W^-}_{j}^{n} \end{bmatrix}\\
&& +  {\Delta t}\frac{1}{\xi} \begin{bmatrix} {\Psi^+}_j^n \\ {\Psi^-}_j^n \end{bmatrix}^\top \begin{bmatrix} \left(I_{m} - {\frac{\Delta t}{\Delta x}}{\Lambda^+}_{j-1}\right) P_j^+ & 0 \\ 0 & \left(I_{k-m} - {\frac{\Delta t}{\Delta x}}{\Lambda^-}_{j+1}\right) P_j^- \end{bmatrix}\begin{bmatrix} {\Psi^+}_j^n \\ {\Psi^-}_j^n \end{bmatrix} \\
&& + {\Delta t} \xi \begin{bmatrix} {W^+}_{j-1}^n \\ {W^-}_{j+1}^n \end{bmatrix}^{\top}  \begin{bmatrix} {\frac{\Delta t}{\Delta x}}{\Lambda_{j-1}^+}P_j^+ & 0 \\ 0 & {\frac{\Delta t}{\Delta x}}{\Lambda_{j+1}^-}P_j^- \end{bmatrix} \begin{bmatrix} {W^+}_{j-1}^n \\ {W^-}_{j+1}^n \end{bmatrix} \\
&& + {\Delta t} \frac{1}{\xi} \begin{bmatrix} {\Psi^+}_j^n \\ {\Psi^-}_j^n \end{bmatrix}^{\top}  \begin{bmatrix} {\frac{\Delta t}{\Delta x}}{\Lambda_{j-1}^+}P_j^+ & 0 \\ 0 & {\frac{\Delta t}{\Delta x}}{\Lambda_{j+1}^-}P_j^- \end{bmatrix}\begin{bmatrix} {\Psi^+}_j^n \\ {\Psi^-}_j^n \end{bmatrix}  + \left({\Delta t}\right)^2 {\Psi_j^n}^{\top}P_j\Psi_j^n,\\
&=&W_j^{n \top} P_j W_j^n - (1+ \xi {\Delta t})\begin{bmatrix} {W^+}_{j}^{n}\\{W^-}_{j}^{n} \end{bmatrix}^{\top}  \begin{bmatrix} {\frac{\Delta t}{\Delta x}}{\Lambda_{j-1}^+}P_j^+ & 0 \\ 0 &  {\frac{\Delta t}{\Delta x}}{\Lambda_{j+1}^-}P_j^- \end{bmatrix} \begin{bmatrix} {W^+}_{j}^{n}\\{W^-}_{j}^{n} \end{bmatrix} \\
&& + (1+ \xi {\Delta t})\begin{bmatrix} {W^+}_{j-1}^{n}\\{W^-}_{j+1}^{n} \end{bmatrix}^{\top}  \begin{bmatrix} {\frac{\Delta t}{\Delta x}}{\Lambda_{j-1}^+}P_j^+ & 0 \\ 0 &  {\frac{\Delta t}{\Delta x}}{\Lambda_{j+1}^-}P_j^- \end{bmatrix} \begin{bmatrix} {W^+}_{j-1}^{n}\\{W^-}_{j+1}^{n} \end{bmatrix}\\
&&\begin{aligned}
&-\begin{bmatrix}  {W^+}_{j}^n - {W^+}_{j-1}^n \\ {W^-}_{j+1}^n - {W^-}_{j}^n \end{bmatrix}^{\top} \begin{bmatrix} {\frac{\Delta t}{\Delta x}}{\Lambda_{j-1}^+}\left(I_{m} - {\frac{\Delta t}{\Delta x}}{\Lambda^+}_{j-1}\right) P_j^+ & 0 \\ 0 &  {\frac{\Delta t}{\Delta x}}{\Lambda_{j+1}^-}\left(I_{k-m} - {\frac{\Delta t}{\Delta x}}{\Lambda^-}_{j+1}\right)P_j^- \end{bmatrix}\\
&\times   \begin{bmatrix}  {W^+}_{j}^n - {W^+}_{j-1}^n \\ {W^-}_{j+1}^n - {W^-}_{j}^n \end{bmatrix}
\end{aligned} \nonumber\\
&& + \xi {\Delta t} W_j^{n \top} P_j W_j^n + {\Delta t}\left(\frac{1}{\xi} + {\Delta t}\right) {\Psi_j^n}^{\top}P_j\Psi_j^n.
\end{eqnarray*}
It can thus be concluded that 
\begin{eqnarray}
\widetilde{W}_j^{n \top} P_j \widetilde{W}_j^{n} &\leq & W_j^{n \top} P_j W_j^n - (1+ \xi {\Delta t})\begin{bmatrix} {W^+}_{j}^{n}\\{W^-}_{j}^{n} \end{bmatrix}^{\top}  \begin{bmatrix} {\frac{\Delta t}{\Delta x}}{\Lambda_{j-1}^+}P_j^+ & 0 \\ 0 &  {\frac{\Delta t}{\Delta x}}{\Lambda_{j+1}^-}P_j^- \end{bmatrix} \begin{bmatrix} {W^+}_{j}^{n}\\{W^-}_{j}^{n} \end{bmatrix} \nonumber\\
&& + (1+ \xi {\Delta t})\begin{bmatrix} {W^+}_{j-1}^{n}\\{W^-}_{j+1}^{n} \end{bmatrix}^{\top}  \begin{bmatrix} {\frac{\Delta t}{\Delta x}}{\Lambda_{j-1}^+}P_j^+ & 0 \\ 0 &  {\frac{\Delta t}{\Delta x}}{\Lambda_{j+1}^-}P_j^- \end{bmatrix} \begin{bmatrix} {W^+}_{j-1}^{n}\\{W^-}_{j+1}^{n} \end{bmatrix}\label{eq:DtDiscLyapunovfun02-002}\\
&& + \xi {\Delta t} W_j^{n \top} P_j W_j^n + {\Delta t}\left(\frac{1}{\xi} + {\Delta t}\right) {\Psi_j^n}^{\top}P_j\Psi_j^n.\nonumber
\end{eqnarray}
Thus, from  inequality \eqref{eq:DtDiscLyapunovfun02-002}, for $ n = 0, \dots, N-1$, equation \eqref{eq:DtDiscLyapunovfun02} is approximated as 
\begin{eqnarray}
\frac{\widetilde{\mathcal{L}}^n-\mathcal{L}^{n}}{\Delta t} &\leq & \xi{\Delta x} \sum_{j=0}^{J-1}W_j^{n \top} P_j W_j^n + \left(\frac{1}{\xi}+ {\Delta t}\right){\Delta x} \sum_{j=0}^{J-1} \Psi_j^{n \top}P_j\Psi_j^n\nonumber\\
&& - \left(1 + \xi{\Delta t} \right) \sum_{j=0}^{J-1} \begin{bmatrix} {W^+}_{j}^{n}\\{W^-}_{j}^{n} \end{bmatrix}^{\top} \begin{bmatrix} {\Lambda_{j-1}^+}P_j^+ & 0 \\ 0 &  {\Lambda_{j+1}^-}P_j^- \end{bmatrix} \begin{bmatrix} {W^+}_{j}^{n}\\{W^-}_{j}^{n} \end{bmatrix} \nonumber \\
&& + \left(1 + \xi{\Delta t} \right) \sum_{j=0}^{J-1} \begin{bmatrix} {W^+}_{j-1}^{n}\\{W^-}_{j+1}^{n} \end{bmatrix}^{\top} \begin{bmatrix} {\Lambda_{j-1}^+}P_j^+ & 0 \\ 0 &  {\Lambda_{j+1}^-}P_j^- \end{bmatrix}\begin{bmatrix} {W^+}_{j-1}^{n}\\{W^-}_{j+1}^{n} \end{bmatrix}. \label{eq:DtDiscLyapunovfun02-01}
\end{eqnarray}
By using $ x_j = x_{j-1} + {\Delta x}, \; j= 0,\dots,J-1 $, boundary conditions \eqref{eq:DiscBCsLHSBLs}, the compatibility conditions \eqref{eq:DiscCCs01LHSBLs} and the assumption in Theorem \ref{thm:ThmChap01-02}, we obtain \cite{Banda2018} 
\begin{eqnarray}
&& \sum_{j=0}^{J-1} \begin{bmatrix} {W^+}_{j-1}^{n}\\{W^-}_{j+1}^{n} \end{bmatrix}^{\top} \begin{bmatrix} {\Lambda_{j-1}^+}P_j^+ & 0 \\ 0 &  {\Lambda_{j+1}^-}P_j^- \end{bmatrix} \begin{bmatrix}  {W^+}_{j-1}^n \\ {W^-}_{j+1}^n \end{bmatrix} \nonumber \\
&& \label{eq:DtDiscLyapunovfun02-02} \\
&& \leq  \sum_{j=0}^{J-1} \begin{bmatrix} {W^+}_{j}^{n}\\{W^-}_{j}^{n} \end{bmatrix}^{\top}  \begin{bmatrix} {\Lambda_{j}^+}P_{j+1}^+ & 0 \\ 0 &  {\Lambda_{j}^-}P_{j-1}^- \end{bmatrix} \begin{bmatrix}  {W^+}_{j}^n \\ {W^-}_{j}^n \end{bmatrix} \nonumber
\end{eqnarray}

for  $n = 0, \dots, N-1, \quad j = 0, \dots, J-1$. We now substitute the inequality \eqref{eq:DtDiscLyapunovfun02-02} in inequality \eqref{eq:DtDiscLyapunovfun02-01} to obtain
\begin{eqnarray}
\frac{\widetilde{\mathcal{L}}^n-\mathcal{L}^{n}}{\Delta t} &\leq &\; \xi{\Delta x} \sum_{j=0}^{J-1}W_j^{n \top} P_j W_j^n + \left(\frac{1}{\xi}+ {\Delta t}\right){\Delta x} \sum_{j=0}^{J-1} \Psi_j^{n \top}P_j\Psi_j^n\nonumber\\
&& - \left(1 + \xi{\Delta t} \right) \sum_{j=0}^{J-1} \begin{bmatrix} {W^+}_{j}^{n}\\{W^-}_{j}^{n} \end{bmatrix}^{\top} \begin{bmatrix} {\Lambda_{j-1}^+}P_j^+ & 0 \\ 0 &  {\Lambda_{j+1}^-}P_j^- \end{bmatrix} \begin{bmatrix} {W^+}_{j}^{n}\\{W^-}_{j}^{n} \end{bmatrix} \nonumber\\
&& + \left(1 + \xi{\Delta t} \right) \sum_{j=0}^{J-1} \begin{bmatrix} {W^+}_{j}^{n}\\{W^-}_{j}^{n} \end{bmatrix}^{\top}  \begin{bmatrix} {\Lambda_{j}^+}P_{j+1}^+ & 0 \\ 0 &  {\Lambda_{j}^-}P_{j-1}^- \end{bmatrix} \begin{bmatrix}  {W^+}_{j}^n \\ {W^-}_{j}^n \end{bmatrix},\nonumber\\
&=& \xi{\Delta x} \sum_{j=0}^{J-1}W_j^{n \top} P_j W_j^n + \left(\frac{1}{\xi}+ {\Delta t}\right){\Delta x} \sum_{j=0}^{J-1} \Psi_j^{n \top}P_j\Psi_j^n\nonumber\\
&& - \left(1 + \xi{\Delta t} \right){\Delta x} \sum_{j=0}^{J-1} \begin{bmatrix} {W^+}_{j}^{n}\\{W^-}_{j}^{n} \end{bmatrix}^{\top} \begin{bmatrix} \frac{{\Lambda_{j-1}^+}P_j^+ - {\Lambda_{j}^+}P_{j+1}^+}{\Delta x} & 0 \\ 0 &  \frac{{\Lambda_{j+1}^-}P_j^- - {\Lambda_{j}^-}P_{j-1}^-}{\Delta x} \end{bmatrix} \begin{bmatrix} {W^+}_{j}^{n}\\{W^-}_{j}^{n} \end{bmatrix},\nonumber\\
&=&\; \xi{\Delta x} \sum_{j=0}^{J-1}W_j^{n \top} P_j W_j^n + \left(\frac{1}{\xi}+ {\Delta t}\right){\Delta x} \sum_{j=0}^{J-1} \Psi_j^{n \top}P_j\Psi_j^n\nonumber\\
&& - \left(1 + \xi{\Delta t} \right){\Delta x} \sum_{j=0}^{J-1} \begin{bmatrix} {W^+}_{j}^{n}\\{W^-}_{j}^{n} \end{bmatrix}^{\top} \begin{bmatrix} -{\Lambda_{j-1}^+}\left(\frac{P_{j+1}^+ - P_{j}^+}{\Delta x} \right)  & 0 \\ 0 &  {\Lambda_{j+1}^-}\left(\frac{P_{j}^- - P_{j-1}^-}{\Delta x} \right)  \end{bmatrix} \begin{bmatrix} {W^+}_{j}^{n}\\{W^-}_{j}^{n} \end{bmatrix}\nonumber\\
&& - \left(1 + \xi{\Delta t} \right){\Delta x} \sum_{j=0}^{J-1} \begin{bmatrix} {W^+}_{j}^{n}\\{W^-}_{j}^{n} \end{bmatrix}^{\top} \begin{bmatrix} -\left(\frac{{\Lambda_{j}^+} - {\Lambda_{j-1}^+}}{\Delta x}\right)P_{j+1}^+ & 0 \\ 0 &   \left(\frac{{\Lambda_{j+1}^-} - {\Lambda_{j}^-}}{\Delta x}\right)P_{j-1}^- \end{bmatrix} \begin{bmatrix} {W^+}_{j}^{n}\\{W^-}_{j}^{n} \end{bmatrix},\nonumber\\
&=& - {\Delta x} \sum_{j=0}^{J-1}W_j^{n \top} \Theta_j W_j^n + \left(\frac{1}{\xi}+ {\Delta t}\right){\Delta x} \sum_{j=0}^{J-1} \Psi_j^{n \top}P_j\Psi_j^n, \label{eq:DtDiscLyapunovfun-002} 
\end{eqnarray}
for $n = 0, \dots, N-1$ where 
\begin{multline*}
\Theta_j := -\left(1 + \xi {\Delta t} \right)\begin{bmatrix} {\Lambda_{j-1}^+}\left(\frac{P_{j+1}^+ - P_{j}^+}{\Delta x} \right)  & 0 \\ 0 &  -{\Lambda_{j+1}^-}\left(\frac{P_{j}^- - P_{j-1}^-}{\Delta x} \right)  \end{bmatrix}\\
-\left(1 + \xi {\Delta t} \right) \begin{bmatrix} \left(\frac{{\Lambda_{j}^+} - {\Lambda_{j-1}^+}}{\Delta x}\right)P_{j+1}^+ & 0 \\ 0 &   -\left(\frac{{\Lambda_{j+1}^-} - {\Lambda_{j}^-}}{\Delta x}\right)P_{j-1}^- \end{bmatrix} -\xi P_j,\; j = 0, \dots, J-1.
\end{multline*}

We use inequality \eqref{eq:DiscLem01-cond02} to obtain  
\begin{align}
\frac{\widetilde{\mathcal{L}}^n-\mathcal{L}^{n}}{\Delta t} \leq &\; - {\Delta x} \sum_{j=0}^{J-1}W_j^{n \top} \Theta_j W_j^n  + {\beta}\left(\frac{1}{\xi}+ {\Delta t}\right) \sup_{0 \leq s \leq n}\left({\Delta x}\sum_{j=0}^{J-1} |\Psi_j^s|^2\right),\label{eq:DtDiscLyapunovfun-0002}
\end{align}
for $n = 0, \dots, N-1$. Assume that $ {\Theta}_j,\ j = 0, \dots, J-1$, is a positive definite matrix. By this assumption, there exist a positive real number $ \eta > 0 $ ($ \eta $ is explicitly defined for specific examples in Section \ref{sec:sec03}) such that $ W_j^{n \top} \Theta_j W_j^n \geq \eta W_j^{n \top} P_j W_j^n $, for $ j = 0, \dots, J-1 $. Therefore, from inequality \eqref{eq:DtDiscLyapunovfun-0002}, inequality \eqref{eq:DiscLem01-cond01} can be obtained as:
\begin{align}
\frac{\widetilde{\mathcal{L}}^n-\mathcal{L}^{n}}{\Delta t}  &\leq - {\eta}\mathcal{L}^{n} + {\beta}\left(\frac{1}{\xi}+ {\Delta t}\right) \sup_{0 \leq s \leq n}\left({\Delta x}\sum_{j=0}^{J-1} |\Psi_j^s|^2\right).\label{eq:DiscLem01-cond01}
\end{align}
Hence, by combining the inequalities \eqref{eq:DtDiscLyapunovfun01} and \eqref{eq:DiscLem01-cond01}, the inequality \eqref{eq:DtDiscLyapunovfun} is approximated as 
\begin{equation}\label{eq:DtDiscLyapunovfun0}
\frac{\mathcal{L}^{n+1} - \mathcal{L}^{n}}{\Delta t} \leq - {\eta}\mathcal{L}^{n} + {\beta}\left(\frac{1}{\xi}+ {\Delta t}\right) \sup_{0 \leq s \leq n}\left({\Delta x}\sum_{j=0}^{J-1} |\Psi_j^s|^2\right),\; n = 0, \dots, N-1.
\end{equation}

By applying Lemma \ref{lem:DiscLemma2} in inequality \eqref{eq:DtDiscLyapunovfun0} and $(1 - \eta \Delta{t})^{n+1} \le e^{-\eta t^{n+1}}$, we have 
\begin{align}
\mathcal{L}^{n+1} \leq &\; \left(\mathcal{L}^{0} - \frac{\beta}{\eta}\left(\frac{1}{\xi}+ {\Delta t}\right) \sup_{0 \leq s \leq n}\left({\Delta x}\sum_{j=0}^{J-1} |\Psi_j^s|^2\right)\right)\left(1 - \eta {\Delta t} \right)^{n+1}\nonumber \\
&\;+ \frac{\beta}{\eta}\left(\frac{1}{\xi}+ {\Delta t}\right) \sup_{0 \leq s \leq n}\left({\Delta x}\sum_{j=0}^{J-1} |\Psi_j^s|^2\right),\nonumber\\
& \leq \; e^{-{\eta}{t}^{n+1}}\mathcal{L}^{0} + \frac{\beta}{\eta}\left(\frac{1}{\xi}+ {\Delta t}\right) \sup_{0 \leq s \leq n}\left({\Delta x}\sum_{j=0}^{J-1} |\Psi_j^s|^2\right), \; n = 0, \dots, N-1.\label{eq:DtDiscLyapunovfun-04}
\end{align}

Thus, from the inequalities \eqref{eq:DtDiscLyapunovfun-04} and \eqref{eq:DiscLem01-cond02}, for all $ j = 0, \dots, J-1 $, $ n = 0, \dots, N-1 $, we get 
\begin{equation}\label{eq:DtDiscLyapunovfun-05}
	\zeta {\Delta x }\sum_{j=0}^{J-1} |W_j^{n+1}|^2 \leq \beta e^{-{\eta}{t}^{n+1}} {\Delta x }\sum_{j=0}^{J-1} |W_j^0|^2 + \frac{\beta}{\eta}\left(\frac{1}{\xi}+ {\Delta t}\right) \sup_{0 \leq s \leq n}\left({\Delta x}\sum_{j=0}^{J-1} |\Psi_j^s|^2\right).  
\end{equation}

From inequality \eqref{eq:DtDiscLyapunovfun-05}, one observes that for $ C = \beta/\zeta $, the condition for the discrete ISS, Equation \eqref{eq:DiscExponStabCondChap1}, is satisfied. Hence, the steady-state $ W_j^n \equiv 0,\; j = 0, \dots, J-1,\; n = 0, \dots, N-1 $ of system \eqref{eq:DiscLHSBLs} with boundary conditions \eqref{eq:DiscBCsLHSBLs} is discrete ISS in $L^2-$norm with respect to discrete disturbance function $ \Psi_j^n,\;j=0,\dots, J-1 $. \qed 
\end{Pf}

The discretisation of the system \eqref{eq:LHSBLs2} is expressed as 
\begin{subequations}\label{eq:DiscLHSBLs2}
\begin{align}
	&\begin{bmatrix} \widetilde{W^+}_{j}^{n}\\\widetilde{W^-}_{j}^{n} \end{bmatrix} =  \begin{bmatrix} {W^+}_{j}^{n}\\{W^-}_{j}^{n} \end{bmatrix}  - {\frac{\Delta t}{\Delta x}} \begin{bmatrix} {\Lambda^+} & 0 \\ 0 & - {\Lambda^-} \end{bmatrix} \begin{bmatrix}  {W^+}_{j}^n - {W^+}_{j-1}^n \\ {W^-}_{j+1}^n - {W^-}_{j}^n \end{bmatrix} + {\Delta t}\begin{bmatrix} {\Psi^+}_j^n \\ {\Psi^-}_j^n \end{bmatrix}, \label{eq:DiscLHSBLs01cc}\\
	&\begin{bmatrix} {W^+}_{j}^{n+1}\\{W^-}_{j}^{n+1} \end{bmatrix} =  
	\begin{bmatrix} \widetilde{W^+}_{j}^{n}\\\widetilde{W^-}_{j}^{n} \end{bmatrix} - {\Delta t}\Pi \begin{bmatrix} \widetilde{W^+}_{j}^{n}\\\widetilde{W^-}_{j}^{n}  \end{bmatrix},\; j = 0, \dots, J-1,\; n = 0, \dots, N-1.\label{eq:DiscLHSBLs02cc}
\end{align}
\end{subequations}

\begin{Cc}\label{col:Col02}
Assume system \eqref{eq:DiscLHSBLs2} with boundary conditions \eqref{eq:DiscBCsLHSBLs} satisfies Assumption \eqref{Asmpt:LHSBLs-Asmpt02} for the discretised system \eqref{eq:DiscLHSBLs2}. Let $ T > 0 $ be fixed and the CFL condition, $ \frac{\Delta t}{\Delta x} \max_{1\leq i \leq k} |\lambda_{i}| \leq 1 $ hold. Let $ \xi $ be any positive real number. Define a positive diagonal matrix, $P_j= \text{diag}\{ P_j^+ , P_j^-\}$, where $ P_j^+$ and $P_j^-$ denote the first $ m $ and the last $ k-m $  diagonal entries, respectively, for $j = 0, \dots, J-1$. Assume that the matrix 
\begin{equation}\label{eq:DiscCol01-cond00}
-\left(1 + \xi {\Delta t} \right)\begin{bmatrix} {\Lambda^+}\left(\frac{P_{j+1}^+ - P_{j}^+}{\Delta x} \right)  & 0 \\ 0 &  -{\Lambda^-}\left(\frac{P_{j}^- - P_{j-1}^-}{\Delta x} \right)  \end{bmatrix} -\xi P_j,
\end{equation}
is positive definite for all $ j = 0, \dots, J-1$ and the matrices
\begin{equation}\label{eq:DiscCol01-cond01}
	P_j\Pi  + {\Pi}^{\top} P_j - {\Delta t}{\Pi_j}^{\top}P_j{\Pi_j},
\end{equation}
and 
\begin{equation}\label{eq:DiscCol01-cond02}
	\begin{bmatrix} {\Lambda^+}P_{J}^+ & 0 \\ 0 &  {\Lambda^-}P_{-1}^- \end{bmatrix} - K^{\top}\begin{bmatrix}  {\Lambda^+}P_{0}^+ & 0 \\ 0 &  {\Lambda^-}P_{J-1}^- \end{bmatrix}K,
\end{equation}
are positive semi-definite for all  $ j = 0, \dots, J-1$. Then the discrete Lyapunov function defined by Equation \eqref{eq:DiscLyapunovfunCh01} is a discrete ISS-Lyapunov function for system \eqref{eq:DiscLHSBLs2} with boundary conditions \eqref{eq:DiscBCsLHSBLs}. Moreover, the steady-state $ W_j^n \equiv 0,\; j = 0, \dots, J-1,\; n = 0, \dots, N-1 $ of system \eqref{eq:DiscLHSBLs2} with boundary conditions \eqref{eq:DiscBCsLHSBLs} is discrete ISS in $L^2-$norm with respect to discrete disturbance function $ \Psi_j^n,\;j=0,\dots, J-1 $.
\end{Cc}
The proof of Corollary \ref{col:Col02} follows from the proof of Theorem \ref{thm:ThmChap01-02} for system \eqref{eq:DiscLHSBLs2}. 

\section{Computational applications and results}\label{sec:sec03}
In this section, numerical tests will be undertaken. The theoretical and numerical results of  Section \ref{sec:sec01} and Section \ref{sec:sec02} will be tested on a linear problem and the Saint-Venant equations.

\subsection{Linear Hyperbolic Systems of Balance Laws}\label{subsec:sec03-1}
We consider a non-uniform $2 \times 2 $ linear hyperbolic system of balance laws:
\begin{equation}\label{eq:2by2LHSBLaws}
\partial_t\begin{bmatrix} w_1\\ w_2 \end{bmatrix} + \begin{bmatrix} \lambda_1(x) & 0  \\ 0 & \lambda_2(x) \end{bmatrix} \partial_x \begin{bmatrix} w_1\\ w_2 \end{bmatrix}
+ \begin{bmatrix} \gamma_{11}(x) & \gamma_{12}(x) \\ \gamma_{21}(x)& \gamma_{22}(x) \end{bmatrix}\begin{bmatrix}w_1\\ w_2\end{bmatrix} = \begin{bmatrix} \psi_1 \\ \psi_2 \end{bmatrix},
\end{equation}
$ x \in [0,l],\; t \in [0, +\infty)$, where $ \lambda_{2}(x) < 0 < \lambda_{1}(x) $. Assume that the system \eqref{eq:2by2LHSBLaws} satisfies \eqref{Asmpt:LHSBLs-Asmpt01}. Set an initial condition as
\begin{equation}\label{eq:2by2LHSBLaws-IC}
\begin{bmatrix} w_1(x,0) \\ w_2(x,0) \end{bmatrix}  = \begin{bmatrix}  f(x)\\ g(x) \end{bmatrix},\; x \in (0,l),
\end{equation}
where $ f $ and $ g $ are smooth functions. Define boundary conditions by
\begin{equation}\label{eq:2by2LHSBLaws-BCs}
\begin{bmatrix} w_1(0,t) \\ w_2(l,t) \end{bmatrix}  = \begin{bmatrix}
0 & k_{12} \\ k_{21} & 0 \end{bmatrix}\begin{bmatrix} w_1(l,t) \\ w_2(0,t) \end{bmatrix},\; t \in [0,+\infty),
\end{equation}
and set compatibility conditions as
\begin{equation}\label{eq:2by2LHSBLaws-CCs}
\begin{bmatrix} w_1(0,0) \\ w_2(l,0) \end{bmatrix}  = \begin{bmatrix}
0 & k_{12} \\ k_{21} & 0 \end{bmatrix}\begin{bmatrix} w_1(l,0) \\ w_2(0,0) \end{bmatrix},
\end{equation}
where  $ k_{12} $ and $ k_{21} $ are constant parameters.

A steady-state solution of the system \eqref{eq:2by2LHSBLaws} can be obtained by solving the following linear system of ordinary differential equations with variable coefficients
\begin{equation}\label{eq:2by2HSBLaws-SteadyState}
\dfrac{d}{dx}\begin{bmatrix} w_1^*(x)\\ w_2^*(x) \end{bmatrix} = \begin{bmatrix} -\frac{\gamma_{11}(x)}{\lambda_1(x)} & -\frac{\gamma_{12}(x)}{\lambda_1(x)} \\ - \frac{\gamma_{21}(x)}{\lambda_2(x)}& -\frac{\gamma_{22}(x)}{\lambda_2(x)} \end{bmatrix}\begin{bmatrix}w_1^*(x)\\ w_2^*(x)\end{bmatrix} + \begin{bmatrix}\frac{\psi_1^*(x)}{\lambda_1(x)}  \\ \frac{\psi_2^*(x)}{\lambda_2(x)} \end{bmatrix}, \; x \in [0,l],
\end{equation}
where the non-uniform steady-states, $ w_1^*(x) $ and $ w_2^*(x) $ may be computed by the Wronskian and Liouville's Formula or by the Lagrange Method.

By following the discussion in \eqref{sec:sec02}, the system \eqref{eq:2by2LHSBLaws} can be split and discretised together with the initial condition \eqref{eq:2by2LHSBLaws-IC}, the boundary conditions \eqref{eq:2by2LHSBLaws-BCs} and the compatibility conditions \eqref{eq:2by2LHSBLaws-CCs} as follows 
\begin{subequations}\label{eq:Disc2by2HSBLaws}
\begin{eqnarray}
	&\begin{bmatrix} \widetilde{w_1}_{j}^{n}\\ \widetilde{w_2}_{j}^{n} \end{bmatrix} =  \begin{bmatrix} {w_1}_{j}^{n}\\{w_2}_{j}^{n} \end{bmatrix}  - {\frac{\Delta t}{\Delta x}} \begin{bmatrix}{\lambda}_{1,j-1} & 0\\0 & {\lambda_{2,j+1}} \end{bmatrix} \begin{bmatrix} {w_1}_{j}^{n}-{w_1}_{j-1}^{n}\\{w_2}_{j+1}^{n} - {w_2}_{j}^{n} \end{bmatrix} + {\Delta t}\begin{bmatrix} {\psi_1}_{j}^{n} \\ {\psi_2}_{j}^{n} \end{bmatrix},\label{eq:Disc2by2HSBLaws01}\\
	&\begin{bmatrix} {w_1}_{j}^{n+1}\\{w_2}_{j}^{n+1} \end{bmatrix} = \begin{bmatrix} \widetilde{w_1}_{j}^{n}\\\widetilde{w_2}_{j}^{n} \end{bmatrix} - {\Delta t} \begin{bmatrix} \gamma_{11,j} & \gamma_{12,j} \\ \gamma_{21,j} & \gamma_{22,j} \end{bmatrix} \begin{bmatrix} \widetilde{w_1}_{j}^{n}\\ \widetilde{w_2}_{j}^{n} \end{bmatrix},\label{eq:Disc2by2HSBLaws02}
\end{eqnarray}
for $ n = 0, \dots , N-1$ and $j = 0, \dots, J-1$,

	\begin{equation}\label{eq:Disc2by2HSBLawsIC}
	w_{1,j}^0 = f_j,\quad w_{2,j}^0 = g_j,\; j = 0, \dots,J-1,
	\end{equation}
	\begin{equation}\label{eq:Disc2by2HSBLawsBCs}
	\begin{bmatrix} {w_1}_{-1}^{n+1}\\{w_2}_{J}^{n+1} \end{bmatrix} =  \begin{bmatrix}0 & k_{12}\\k_{21} & 0 \end{bmatrix} \begin{bmatrix} {w_1}_{J-1}^{n+1}\\{w_2}_{0}^{n+1} \end{bmatrix},\quad n = 0, \dots , N-1,
	\end{equation}
	\begin{equation}\label{eq:Disc2by2HSBLawsCCs}
	\begin{bmatrix} {w_1}_{-1}^{0}\\{w_2}_{J}^{0} \end{bmatrix} =  \begin{bmatrix}0 & k_{12}\\k_{21} & 0 \end{bmatrix} \begin{bmatrix} {w_1}_{J-1}^{0}\\{w_2}_{0}^{0} \end{bmatrix}.
	\end{equation}
\end{subequations}

For a fixed $ T > 0 $, we apply the CFL condition: \[\frac{\Delta t}{\Delta x}\max_{0 \leq j \leq J-1}\{|\lambda_{1,j}|,|\lambda_{2,j}|\} \leq 1. \]  
Now apply the $ L^2- $ ISS-Lyapunov function \eqref{eq:DiscLyapunovfunCh01} for system \eqref{eq:Disc2by2HSBLaws} and consider the assumptions of Theorem \ref{thm:ThmChap01-02}:
\begin{itemize}
	\item[\textbf{C1:}] the matrix 
	\begin{multline*}
	{\theta}_{j} := -\left(1 + \xi {\Delta t} \right)\begin{bmatrix} {{\lambda_1}_{j-1}}\left(\frac{{p_1}_{j+1} - {p_1}_{j}}{\Delta x} \right)  & 0 \\ 0 &  {{\lambda_2}_{j+1}}\left(\frac{{p_2}_{j} - {p_2}_{j-1}}{\Delta x} \right)  \end{bmatrix}\\
	-\left(1 + \xi {\Delta t} \right) \begin{bmatrix} \left(\frac{{{\lambda_1}_{j}} - {{\lambda_1}_{j-1}}}{\Delta x}\right){p_1}_{j+1} & 0 \\ 0 &   \left(\frac{{{\lambda_2}_{j+1}} - {{\lambda_2}_{j}}}{\Delta x}\right){p_2}_{j-1} \end{bmatrix} -\xi P_j,
	\end{multline*}
	is positive definite for all $ j = 0, \dots, J-1$,
	\item[\textbf{C2:}] the matrix 
	\begin{multline*}
	M_{j} := \begin{bmatrix} {p_1}_{j} & 0 \\ 0 & {p_2}_{j} \end{bmatrix}\begin{bmatrix} {\gamma_{11}}_{j} & {\gamma_{12}}_{j} \\ {\gamma_{21}}_{j} & {\gamma_{22}}_{j} \end{bmatrix} + \begin{bmatrix} {\gamma_{11}}_{j} & {\gamma_{12}}_{j} \\ {\gamma_{21}}_{j} & {\gamma_{22}}_{j} \end{bmatrix}^{\top} \begin{bmatrix} {p_1}_{j} & 0 \\ 0 & {p_2}_{j} \end{bmatrix}\nonumber\\
	-{\Delta t} \begin{bmatrix} {\gamma_{11}}_{j} & {\gamma_{12}}_{j} \\ {\gamma_{21}}_{j} & {\gamma_{22}}_{j} \end{bmatrix}^{\top}\begin{bmatrix} {p_1}_{j} & 0 \\ 0 & {p_2}_{j} \end{bmatrix}\begin{bmatrix} {\gamma_{11}}_{j} & {\gamma_{12}}_{j} \\ {\gamma_{21}}_{j} & {\gamma_{22}}_{j} \end{bmatrix},
	\end{multline*} 
	is positive semi-definite for all $ j = 0, \dots, J-1$, and
	\item[\textbf{C3:}] the matrix 
	\begin{eqnarray*}
	B_c:=&& \begin{bmatrix} {\lambda_1}_{J-1} {p_1}_{J} & 0 \\ 0 & |{\lambda_2}_{0}| {p_2}_{-1}\end{bmatrix}	\\
	&& - \begin{bmatrix}0 & k_{12}\\k_{21} & 0 \end{bmatrix}^{\top} \begin{bmatrix} {\lambda_1}_{-1} {p_1}_{0} & 0 \\ 0 & |{\lambda_2}_{J}| {p_2}_{J-1}\end{bmatrix} \begin{bmatrix}0 & k_{12}\\k_{21} & 0 \end{bmatrix},
	\end{eqnarray*}
	is positive semi-definite.
\end{itemize}

Now we verify the above assumptions. For assumption \textbf{C1} it suffices to show that both diagonal entries of $ \theta_{j} $ are positive, i.e., for all $ j = 0, \dots, J-1$,
\begin{eqnarray*}
&{\eta_1}_{j} := \left(-\left(1 + \xi {\Delta t} \right)\frac{{\lambda_1}_{j-1}}{{p_1}_{j}}\left(\frac{{p_1}_{j+1} - {p_1}_{j}}{\Delta x} \right) -\left(1 + \xi {\Delta t} \right) \left(\frac{{{\lambda_1}_{j}} - {{\lambda_1}_{j-1}}}{\Delta x}\right)\frac{{p_1}_{j+1}}{{p_1}_{j}} -\xi\right) {p_1}_{j} > 0, \\
&{\eta_2}_{j}:= \left(-\left(1 + \xi {\Delta t} \right)  \frac{{\lambda_2}_{j+1}}{{p_2}_{j}}\left(\frac{{p_2}_{j} - {p_2}_{j-1}}{\Delta x} \right) -\left(1 + \xi {\Delta t} \right)  \left(\frac{{{\lambda_2}_{j+1}} - {{\lambda_2}_{j}}}{\Delta x}\right)\frac{{p_2}_{j-1}}{{p_2}_{j}}  -\xi\right) {p_2}_{j} > 0.
\end{eqnarray*}

The second assumption \textbf{C2}  holds if the matrix $ M_{j} $, 
\begin{equation*}
M_{j} = \begin{bmatrix} {M_{11}}_{j} & {M_{12}}_{j} \\ {M_{12}}_{j} & {M_{22}}_{j} \end{bmatrix},\;j = 0, \dots, J-1, 
\end{equation*}
where 
\begin{align*}
M_{11,j} &:= 2 {\gamma_{11}}_{j}{p_1}_{j}-{\Delta t} \left( {\gamma_{11}^{2}}_{j}{p_1}_{j} + {\gamma_{21}^{2}}_{j} {p_2}_{j}\right), \\
M_{12,j}&:=  {\gamma_{21}}_{j}{p_2}_{j}+{\gamma_{12}}_{j}{p_1}_{j}
- {\Delta t} \left( {\gamma_{11}}_{j} {\gamma_{12}}_{j}{p_1}_{j} + {\gamma_{21}}_{j}{\gamma_{22}}_{j}{p_2}_{j} \right),\\
M_{22,j}&:= 2 {\gamma_{22}}_{j}{p_2}_{j} - {\Delta t} \left( {\gamma_{12}^{2}}_{j}{p_1}_{j} + {\gamma_{22}^{2}}_{j}{p_2}_{j}\right),  
\end{align*} 
has non-negative eigenvalues, 
\begin{equation*}
\sigma_j^\pm = \frac{1}{2}\left( \left({M_{11}}_{j} + {M_{22}}_{j}\right) \pm \sqrt{\left({M_{11}}_{j} + {M_{22}}_{j}\right)^2 - 4\left({M_{11}}_{j} {M_{22}}_{j} - {M_{12}^2}_{j}\right)}\right) \geq 0,
\end{equation*}
for all  $ j = 0, \dots, J-1$.  

Finally, the matrix $ B_c $ can be expressed as 
\begin{equation*}
B_c = \begin{bmatrix} {\lambda_1}_{J-1} {p_1}_{J} - k_{21}^2 |{\lambda_2}_{J}| {p_2}_{J-1} & 0 \\ 0 & |{\lambda_2}_{0}| {p_2}_{-1} - k_{12}^2 {\lambda_1}_{-1} {p_1}_{0} \end{bmatrix},
\end{equation*}
and if we can choose the parameters, $ \kappa_{12} $ and $ \kappa_{21} $ as 
\begin{equation*}
k_{12}^2  \leq \frac{|{\lambda_2}_{0}| {p_2}_{-1}}{{\lambda_1}_{-1} {p_1}_{0}}, \quad \text{and} \quad k_{21}^2  \leq \frac{{\lambda_1}_{J-1} {p_1}_{J}}{|{\lambda_2}_{J}| {p_2}_{J-1}},
\end{equation*}
then the assumption \textbf{C3} holds. 

Thus, the approximation of the time derivative of the candidate ISS-Lyapunov function defined by Equation \eqref{eq:DtDiscLyapunovfun} can be expressed as in Equation \eqref{eq:DtDiscLyapunovfun0} with 
$\displaystyle \eta := \min_{0 \leq j \leq J-1} \{{\eta_1}_{j}, {\eta_2}_{j} \} $ and $\displaystyle \beta  = \max_{0 \leq j \leq J-1}P_j $. Hence the candidate ISS-Lyapunov function satisfies Definition \ref{def:DiscISSLyafun} and defines an upper bound for the discrete ISS-Lyapunov function by \eqref{eq:DtDiscLyapunovfun-04}.

We now analyse ISS for the following  linear hyperbolic system of balance laws with additive disturbance in one space dimension \eqref{eq:2by2LHSBLaws} with $\lambda_1(x) = 1$, $\lambda_1(x) = -1$, $\psi_1 = \psi$, $\psi_2 = -\psi$,
where $ \psi $ is defined by 
\begin{equation*}
\psi(x,t) = \begin{dcases} 0.01 \sin^2(\pi t), & 0 \leq t < 5,\\
0, &  t \geq 5.
\end{dcases}
\end{equation*} 

The initial conditions \eqref{eq:2by2LHSBLaws-IC} are set using $f(x) = -0.5$ and $g(x) = 0.5$ for all $x \in [0, 1]$. In addition, boundary conditions and compatibility conditions are defined by \eqref{eq:2by2LHSBLaws-BCs} and \eqref{eq:2by2LHSBLaws-CCs}, respectively,  for $l = 1$. 
Similarly, the discretised system \eqref{eq:Disc2by2HSBLaws} will be considered for the uniform system.

Let the CFL condition, $ \lambda\frac{\Delta t}{\Delta x} \leq 1 $, where $ \lambda = \max \{ \lambda_1, |\lambda_2|\} = 1$ holds for a fixed $ T > 0 $. Define an implicit discrete weight function by \[ \displaystyle P_j := \text{diag}\{ {p_1}\exp(-\mu x_{j}) , {p_2}\exp(\mu x_{j}) \}, \; p_1 > 0, p_2 > 0,  \mu > 0,\; j = 0, \dots, J-1.\] Then, for $ \xi > 0 $, we can choose sufficiently small $ \mu > \xi$ such that for $ p_1 = p_2 = 1 $, we have  $|k_{12}| \leq  1$ and $|k_{21}| \leq \exp(-\mu)$. Therefore, the conditions in Equation \eqref{col:Col02} are satisfied.  Hence, the discrete system with requisite initial conditions, boundary conditions and compatibility conditions of the considered example is discrete ISS for the discrete  $ L^2-$norm. 

For numerical computations, we take CFL = 0.75, $ {\Delta x} = 1/1600 $, $ {\Delta t} = 0.75/1600 $ and $ \xi = 0.125 $. Then, the decay rate 
\begin{eqnarray*}
\eta &=&\; \min\{\eta_1 , \eta_2 \},\\  
	 &=&\; \min \left\lbrace \left(1 + \xi {\Delta t}\right) \lambda_{1} \left(\frac{1 - \exp\left(-\mu {\Delta x}\right)}{\Delta x}\right) - \xi, \left(1 + \xi {\Delta t}\right) |\lambda_{2}| \left(\frac{1 - \exp\left(-\mu {\Delta x}\right)}{\Delta x}\right) - \xi \right\rbrace,\\
	 &=&\; \alpha \left(1 + \xi {\Delta t}\right)\left(\frac{1 - \exp\left(-\mu {\Delta x}\right)}{\Delta x}\right)  - \xi,\\
	 &\geq&\;  \mu \alpha \left(1 + \xi {\Delta t}\right)\exp\left(-\mu {\Delta x}\right) - \xi > 0, \quad \text{if} \quad \xi < \mu < 19098.926,
\end{eqnarray*}
where $\alpha = \min \{\lambda_{1}, |\lambda_{2}| \} = 1$. Here, $|k_{21}| \leq \exp(-\mu) \approx 0.8825$.

For the above choice of values and if $ p_1 = p_2 = 1 $ is chosen, the assumptions \textbf{C1-C3} hold. Furthermore, the upper bound of the discrete ISS-Lyapunov function is defined by Equation \eqref{eq:DtDiscLyapunovfun-04} with 
\[
\beta  = \max\left\{ \max_{ 0 \leq j \leq J-1}\{p_1 \exp(-\mu x_{j})\},  \max_{ 0 \leq j \leq J-1}\{p_2 \exp(\mu x_{j})\} \right\}.
\]

In addition, we compute a comparison of the discrete ISS-Lyapunov function and its upper bound for CFL = 0.75 and CFL = 1 in Table \ref{table:Convergence-01} and Table \ref{table:Convergence-02}, respectively. 
\begin{table}[H]
	\centering
	\begin{tabular}{ccccccc}
		\hline 
		$ J $ & $ \|\mathcal{L}_{\text{up}}^n - \mathcal{L}^n \|_{L^\infty} $ & $ \|\mathcal{L}_{\text{up}}^n - \mathcal{L}^n \|_{L^2} $ &  $\mu$ & $\eta $ \\[0.5ex] 
		\hline 
		200    & 0.27569    & 0.46295    & 0.575    & 0.44862 \\[1ex]
		400    & 0.27362    & 0.46051    & 0.575    & 0.44931 \\[1ex]
		800    & 0.27218    & 0.45877    & 0.575    & 0.44965 \\[1ex]
		1600   & 0.27119    & 0.45753    & 0.575    & 0.44983 \\[0.5ex]
		\hline 
	\end{tabular} 
	\caption[]%
	{The comparison of the upper bound of the Lyapunov function with discrete Lyapunov function.  Under CFL = 0.75, $ {\Delta x} = \frac1J $, $ {\Delta t} =  \frac{{\Delta x}}{\max \{\lambda_{1}, |\lambda_{2}| \}}\text{CFL}$, $  \xi = 0.125 $,  $ T = 10 $ and $ k_{12} = k_{21} = 0.5$.}
	\label{table:Convergence-01}
\end{table}
 
\begin{table}[H]
	\centering
	\begin{tabular}{ccccccc}
		\hline 
		$ J $ & $ \|\mathcal{L}_{\text{up}}^n - \mathcal{L}^n \|_{L^\infty} $ & $ \|\mathcal{L}_{\text{up}}^n - \mathcal{L}^n \|_{L^2} $ &  $\mu$ & $\eta $ \\[0.5ex] 
		\hline 
		200     & 0.269      & 0.39341    & 0.575    & 0.44871 \\[1ex]
		400     & 0.26891    & 0.39352    & 0.575    & 0.44935 \\[1ex]
		800     & 0.26887    & 0.39357    & 0.575    & 0.44968 \\[1ex]
		1600    & 0.26885    & 0.3936     & 0.575    & 0.44984 \\[0.5ex]
		\hline 
	\end{tabular} 
	\caption[]%
	{The comparison of the upper bound of the Lyapunov function with discrete Lyapunov function.  Under CFL = 1, $ {\Delta x} = \frac1J $, $ {\Delta t} =  \frac{{\Delta x}}{\max \{\lambda_{1}, |\lambda_{2}| \}}\text{CFL}$, $  \xi = 0.125 $,  $ T = 10 $ and $ k_{12} = k_{21} = 0.5$.}
	\label{table:Convergence-02}
\end{table} 
From the results listed in the two tables above, it can be verified that the required estimates are obtained.

\subsection{The Saint-Venant Equations}
During a rainfall or an evaporation, a flow of water along a channel can be affected by an inflow or outflow of water into the channel which changes the depth and velocity of water. As a result, the flow will be disturbed. The dynamics of a water flow along a pool of prismatic horizontal open channel with a rectangular cross section, a unit width, a constant bottom slope $ S_b $ with disturbance (rain is considered) is described by Saint-Venant equations obtained from \cite{bastin2008using, kirstetter2016modeling} as 
\begin{equation}\label{eq:SVEqns}
\begin{split}
	&\partial_t h + \partial_x\left(hu\right) = R,\\
	&\partial_t u + \partial_x\left(\frac{1}{2}u^2 + gh \right) + g\left( C_f\frac{u^2}{h} -  S_b\right) = -\frac{u}{h}R, \; x \in [0,l], \; t \in [0, +\infty),
\end{split}
\end{equation}
where $ h:= h(x,t) $ is water depth, $ u:= u(x,t) $ is water velocity, $ g $ is gravitational acceleration, $ C_f $ is a friction parameter and $R:=R(x,t) $ is rainfall intensity. We set an initial condition as 
\begin{equation}\label{eq:ICSVEqns}
h(x,0)  = h_0(x), \;   u(x,0) = u_0(x), \; x \in (0,l),
\end{equation}
where $ h_0 $ and $ u_0 $ are smooth functions. We define a linear boundary condition together with compatibility conditions by 
\begin{equation}\label{eq:BCsSVEqns}
u(0,t) = \kappa_0 h(0,t), \; u(l,t) = \kappa_lh(l,t), \; t \in [0, +\infty),
\end{equation}
where $ \kappa_0 $ and $ \kappa_l $ are constant parameters.

For any smooth solution of the Saint-Venant equations \eqref{eq:SVEqns}, a spatially dependent steady-state $ h^*(x)$, $ u^*(x) $ satisfies 
\begin{equation}\label{eq:SVEqns-SS01}
\begin{split}
	&\left(h^*(x)u^*(x)\right)' = R^*(x),\\
	&\left(\frac{1}{2}{u^*}^2(x) + gh^*(x) \right)' + g\left( C_f\frac{{u^*}^2(x)}{h^*(x)} -  S_b\right) = -\frac{u^*(x)}{h^*(x)}R^*(x), \; x \in [0,l].
\end{split}
\end{equation}
By solving the first order system of ODEs \eqref{eq:SVEqns-SS01}, we obtain 
\begin{align}\label{eq:SVEqns-SS02}
{h^*}'(x) =&\; \frac{1}{{u^*}^2(x) - gh^*(x)} \left[ gh^*(x)\left(C_f\frac{{u^*}^2(x)}{h^*(x)} -  S_b\right) + 2u^*(x)R^*(x)\right],\\
{u^*}'(x) =&\; \frac{-1}{{u^*}^2(x) - gh^*(x)} \left[ gu^*(x)\left(C_f\frac{{u^*}^2(x)}{h^*(x)} -  S_b\right) + \left(g + \frac{{u^*}^2(x)}{h^*(x)} \right)R^*(x)\right], 
\end{align}
where a sub-critical flow is assumed, i.e. $ gh^*(x) - {u^*}^2(x) > 0,\; x \in [0,l]$. 

For the sub-critical flow, the system \eqref{eq:SVEqns} is strictly hyperbolic since the Jacobian matrix of the flux function has two distinct eigenvalues. Therefore, the system \eqref{eq:SVEqns} can be linearised around the steady-state as follows 
\begin{multline}\label{eq:SVEqns-Linearized}
\begin{bmatrix} v_1\\ v_2 \end{bmatrix}_t + \begin{bmatrix}
	u^*(x) & h^*(x) \\ g & u^*(x) \end{bmatrix}\begin{bmatrix} v_1\\ v_2 \end{bmatrix}_x\\
+ \begin{bmatrix} {u^*}'(x) & {h^*}'(x) \\ -gC_f\frac{{u^*}^2(x)}{{h^*}^2(x)} -\frac{u^*(x)}{{h^*}^2(x)}R^*(x) & {u^*}'(x) + 2gC_f\frac{{u^*}(x)}{{h^*}(x)} + \frac{R^*(x)}{h^*(x)} \end{bmatrix}\begin{bmatrix} v_1\\ v_2 \end{bmatrix} = \begin{bmatrix} \delta \\ -\frac{u^*(x)}{h^*(x)}\delta \end{bmatrix},
\end{multline}
where $ v_1 = h - h^*(x) $, $ v_2 = u - u^*(x) $, $ \delta = R - R^*(x) $ and the Jacobian matrix is diagonalised as 
\begin{equation}\label{eq:SVEqns-Diagoonalized}
H \begin{bmatrix} u^*(x) & h^*(x) \\ g & u^*(x) \end{bmatrix} H^{-1} = \begin{bmatrix} u^*(x) + \sqrt{gh^*(x)} & 0 \\ 0 & u^*(x) - \sqrt{gh^*(x)} \end{bmatrix}.
\end{equation}
where \[ H = \begin{bmatrix}\sqrt{\frac{g}{h^*(x)}} & 1\\ - \sqrt{\frac{g}{h^*(x)}} & 1 \end{bmatrix}.\]
We now define Riemann-invariants (change of coordinates) for the linearised system \eqref{eq:SVEqns-Linearized} by using the diagonalisation \eqref{eq:SVEqns-Diagoonalized} as 
\begin{equation}\label{eq:SVEqns-RiemannIvs01}
\begin{bmatrix} w_1 \\ w_2 \end{bmatrix} = H \begin{bmatrix} v_1\\ v_2 \end{bmatrix} = \begin{bmatrix} v_2 + v_1\sqrt{\frac{g}{H^*(x)}}\\ v_2 - v_1\sqrt{\frac{g}{H^*(x)}} \end{bmatrix}.
\end{equation}
Thus,  
\begin{equation}\label{eq:SVEqns-RiemannIvs02}
\begin{bmatrix} v_1\\ v_2 \end{bmatrix} = H^{-1}\begin{bmatrix} w_1 \\ w_2 \end{bmatrix} = 
\begin{bmatrix} \frac{1}{2}\sqrt{\frac{h^*(x)}{g}}\left( w_1 - w_2\right) \\  \frac{1}{2}\left( w_1 + w_2\right) \end{bmatrix}.
\end{equation}

By using the coordinates expressed in Equation \eqref{eq:SVEqns-RiemannIvs01} or Equation \eqref{eq:SVEqns-RiemannIvs02}, the linearised system \eqref{eq:SVEqns-Linearized} can be decoupled as Equation \eqref{eq:2by2LHSBLaws}
 
where $ \lambda_1(x) = u^*(x) + \sqrt{gh^*(x)} $, $ \lambda_2(x) = u^*(x) - \sqrt{gh^*(x)} $,
\begin{align*}
\gamma_{11}(x)= &\; {u^*}'(x) + \frac{1}{4h^*(x)} \left(\lambda_1(x) + 2\sqrt{gh^*(x)}  \right){h^*}'(x) \\
&\;+ \frac{gC_f{u^*}^2(x)}{2h^*(x)}\left(\frac{2}{u^*(x)} - \frac{1}{\sqrt{gh^*(x)}}\right) -\frac{1}{2h^*(x)\sqrt{gh^*(x)}}\lambda_2(x)R^*(x),\\
\gamma_{12}(x)= &\; -\frac{1}{4h^*(x)} \left(\lambda_1(x) - 2\sqrt{gh^*(x)}  \right){h^*}'(x) \\
&\;+ \frac{gC_f{u^*}^2(x)}{2h^*(x)}\left(\frac{2}{u^*(x)} + \frac{1}{\sqrt{gh^*(x)}}\right) + \frac{1}{2h^*(x)\sqrt{gh^*(x)}}\lambda_1(x)R^*(x),\\
\gamma_{21}(x)= &\; -\frac{1}{4h^*(x)} \left(\lambda_2(x) - 2\sqrt{gh^*(x)}  \right){h^*}'(x) \\
&\;+ \frac{gC_f{u^*}^2(x)}{2h^*(x)}\left(\frac{2}{u^*(x)} - \frac{1}{\sqrt{gh^*(x)}}\right) - \frac{1}{2h^*(x)\sqrt{gh^*(x)}}\lambda_2(x)R^*(x),\\
\gamma_{22}(x)= &\; {u^*}'(x) + \frac{1}{4h^*(x)} \left(\lambda_2(x) - 2\sqrt{gh^*(x)}  \right){h^*}'(x) \\
&\;+ \frac{gC_f{u^*}^2(x)}{2h^*(x)}\left(\frac{2}{u^*(x)} + \frac{1}{\sqrt{gh^*(x)}}\right) + \frac{1}{2h^*(x)\sqrt{gh^*(x)}}\lambda_1(x)R^*(x),\\
\psi_1 = &\;  -\frac{\lambda_2(x)}{h^*(x)} \delta, \quad \psi_2 = \; -\frac{\lambda_1(x)}{h^*(x)} \delta.
\end{align*}
Consequently, the initial condition \eqref{eq:ICSVEqns} and the boundary conditions together with compatibility conditions \eqref{eq:BCsSVEqns} are re-written as  \eqref{eq:2by2LHSBLaws-IC}, \eqref{eq:2by2LHSBLaws-BCs}, \eqref{eq:2by2LHSBLaws-CCs}, respectively, with $f(x) = v_2(x,0) + v_1(x,0)\sqrt{\frac{g}{H^*(x)}}$ and $g(x) = v_2(x,0) - v_1(x,0)\sqrt{\frac{g}{H^*(x)}}$ for $x \in (0,l)$,
\begin{equation*}
k_{12} = \frac{\kappa_0\sqrt{\frac{H^*(0)}{g}} + 1}{\kappa_0\sqrt{\frac{H^*(0)}{g}} - 1}\qquad  \text{ and } \quad k_{21} = \frac{\kappa_1\sqrt{\frac{H^*(l)}{g}} + 1}{\kappa_1\sqrt{\frac{H^*(l)}{g}} - 1}.
\end{equation*}
The linearised and decoupled system is discretised as in Equation \eqref{eq:Disc2by2HSBLaws}.

As a test example, we take a constant steady-state from \cite{diagne2017backstepping}, $ h^*(x) = 2 $, $ u^*(x) = 3 $ for all $ x \in [0,1]$ with 
physical parameters $ g = 9.81 $, $ C_f = 0.1 $ and $ S_b = 0.0459 $. Then, by solving  the steady-state system \eqref{eq:SVEqns-SS01}, we obtain $ R^*(x) = 0 $ for all $x \in [0,1]$. We take a homogeneous rainfall intensity as
\begin{equation*}
R(x,t) = \begin{dcases} 0.25 \sin^2(\pi t), & 0 \leq t < 5,\\
	0, &  t \geq 5.
\end{dcases}
\end{equation*} 
The initial condition for the system \eqref{eq:SVEqns} is taken as: \[H(x,0) = 2.5, \quad V(x,0) = 4\sin(\pi x), \text{for }  x \in [0,1].\]

The linear system has eigenvalues, $ \lambda_1 = 7.4294$ and $\lambda_2 = -1.4294$ and coefficients of the source terms are 
$\gamma_{11} = \gamma_{21} = 0.0992$ and $\gamma_{12} = \gamma_{22} = 0.2008.$

The initial condition in terms of the new coordinates is $$ w_1(x,0) = -1.8926 + 4\sin(\pi x), \quad w_2(x,0) = -4.1074 + 4\sin(\pi x), $$ for $x \in [0,1]$. 

We take $\text{CFL} = 0.75$, $\Delta x = 1/1600$ and $ \xi = 1/8 $, and analyse the numerical boundary feedback ISS for the implicit discrete weight function defined in Section \ref{subsec:sec03-1}. Then, the decay rate  
\begin{equation*}
\eta = \mu \alpha \left(1 + \xi {\Delta t}\right) \exp\left(-\mu {\Delta x} \right) - \xi > 0,\quad \text{if}\quad 0.087446 < \mu < 2008.457445,
\end{equation*}
where $\alpha = 1.4294$.  

\begin{figure}[H]
\centering
\includegraphics[scale=0.7]{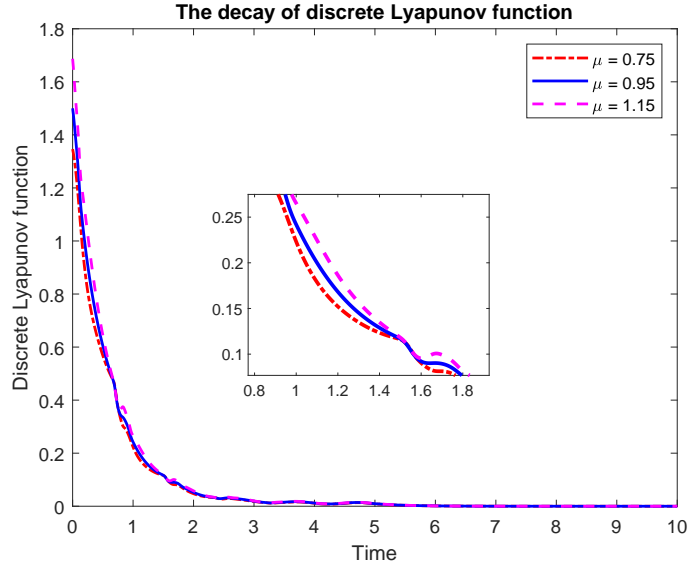}
\caption{The decay of Lyapunov function for Saint-Venant equations. The choice of parameters are $ p_1 =0.0992,  p_2 = 0.2008 $, $ k_{12} = 0.75 $ and $ k_{21} = 0.75 $ with  $ l = 1 $, $ J = 1600 $ and $ T = 10 $ under CFL = 0.75.}\label{fig:SWEwithR}
\end{figure}

For numerical implementation, a sufficiently small value of $\mu$ is chosen such that for the constant steady state the parameters $p_1$ and $p_2$ are chosen to satisfy $p_1\gamma_{12} = p_2\gamma_{21}$. For this example, the values $p_1 = \gamma_{21} = 0.0992$ and $p_2 = \gamma_{12} = 0.2008$ were used. With this choice of parameters, we obtain $|k_{12}| < 0.6241$ and $|k_{21}| < 1.6024e^{-\mu} < 1.4683$. The numerical convergence of the discrete ISS-Lyapunov function for different values of $ \mu $ is shown in Figure \ref{fig:SWEwithR}.

Figure \ref{fig:SWEwithR} illustrates the decay of the ISS-Lyapunov function in the presence of additive disturbance. 

\begin{Rk}
In \cite{bastin2017quadratic}, the explicit weight function, $ P(x) := \frac{1}{2} h^*(x)I_2 $, for Saint-Venant Equations is considered. However, the assumption \textbf{C1} is not satisfied for constant steady-state. Therefore, we may only use it for non-uniform steady-state
\end{Rk}

\section{Conclusion}\label{sec:sec04}

A non-uniform linear hyperbolic system of balance laws with additive disturbance has been considered. A first-order finite volume method with a time-splitting technique is used in the discretisation of this linear system. A theoretical and numerical analysis for boundary control has been presented. An $ \text{L}^2-$ISS-Lyapunov function is used to investigate conditions for ISS of both the continuous and the discretised system.  The decay of the Lyapunov function has been proved. The result was applied to a linear hyperbolic system of balance laws and to a relevant physical problem: the Saint-Venant equations. Explicit computations of the decay have been undertaken and demonstrated and agree with the analytical results. The properties that have been proved using analysis can also be observed in these results.

This result can be used to extend the theory to prove the decay of appropriate Lyapunov functions for non-uniform linear balance laws with boundary disturbance. Such analysis is underway. For a system of the form:
\[\partial_t W + \Lambda (x) \partial_x W + \Pi(x) W = 0,\; x \in [0,l],\; t \in [0,+\infty)\]
with boundary disturbance, preliminary results show that for non-uniform equilibria it is possible to obtain decay of the Lyapunov function in the $L^2$-norm.

This work still leaves more questions open. The problem of analysing ISS-Lyapunov functions for nonlinear hyperbolic differential equations is considered for future work. Further it will also be of interest to consider more accurate finite volume methods. The approach used currently has significant numerical viscosity and this might influence the rate of convergence of the discrete Lyapunov function. Careful analysis of the influence of such numerical artefacts needs to be undertaken.

\section*{Acknowledgments}
This work is supported in part by the National Research Foundation of South Africa (Grant number: 93099 and 102563) and the German Research Foundation (DFG) grant number: GO 1920/10-1.
\bibliographystyle{abbrv}


\end{document}